\documentclass[11pt]{amsart}
\usepackage{amssymb}
\usepackage{latexsym,bm}
\usepackage{amsmath,amsfonts,amssymb,mathrsfs}

\DeclareMathOperator{\End}{End} 
 
 \DeclareMathOperator{\ch}{ch}
\DeclareMathOperator{\cop}{cop} \DeclareMathOperator{\op}{op}
 \newcommand{\lam}{\lambda}
\newcommand{\ksl}{\mathfrak{sl}} \newcommand{\ru}{\rm U}
\newcommand{\DF}{\displaystyle}
\newcommand{\DDF}[2]{\frac{\DF#1}{\DF#2}}
\newcommand{\wtu}{{\rm D}}
\newcommand{\bwtu}{{\mathbb D}}

\newtheorem{proposition}{Proposition}[section]
\newtheorem{lemma}[proposition]{Lemma}
\newtheorem{corollary}[proposition]{Corollary}
\newtheorem{theorem}[proposition]{Theorem}
\newtheorem{remark}[proposition]{Remark}

\theoremstyle{definition}

\newcommand{\thlabel}[1]{\label{th:#1}}
\newcommand{\thref}[1]{Theorem~\ref{th:#1}}

\newcommand{\lelabel}[1]{\label{le:#1}}
\newcommand{\leref}[1]{Lemma~\ref{le:#1}}

\newcommand{\colabel}[1]{\label{co:#1}}
\newcommand{\coref}[1]{Corollary~\ref{co:#1}}

\date{}

\begin{document}
\title{Quantum double of ${\rm U}_q((\ksl_2)^{\leq 0})$}
\author{Jun Hu}
\address{Department of Applied Mathematics\\
Beijing Institute of Technology\\
Beijing, 100081, P.R. China}
\email{junhu303@yahoo.com.cn}
\author{Yinhuo Zhang}
\address{School of Mathematics, Statistics and Computer Science,
Victoria  University of Wellington, PO Box 600, Wellington, New
Zealand} \email{yinhuo.zhang@vuw.ac.nz} 

\begin{abstract} Let ${\rm U}_q(\ksl_2)$ be the quantized enveloping algebra associated to the simple Lie algebra
$\ksl_2$. In this paper, we study the quantum double $\wtu_q$ of the
Borel subalgebra ${\rm U}_q((\ksl_2)^{\leq 0})$ of ${\rm
U}_q(\ksl_2)$. We construct an analogue of Kostant--Lusztig
$\mathbb{Z}[v,v^{-1}]$-form for $\wtu_q$ and show that it is a Hopf
subalgebra. We prove that, over an algebraically closed field, every
simple $\wtu_q$-module is the pullback of a simple ${\rm
U}_q(\ksl_2)$-module through certain surjection from $\wtu_q$ onto
${\rm U}_q(\ksl_2)$, and the category of finite dimensional weight
$\wtu_q$-modules is equivalent to a direct sum of $|k^{\times}|$
copies of the category of finite dimensional weight ${\rm
U}_q(\ksl_2)$-modules. As an application, we recover (in a
conceptual way) Chen's results \cite{C1} as well as Radford's
results \cite{R3} on the quantum double of Taft algebra. Our main
results allow a direct generalization to the quantum double of the
Borel subalgebra of the quantized enveloping algebra associated to
arbitrary Cartan matrix.
\end{abstract}
\maketitle

\section{Preliminaries}

Let $k$ be a field. Let $q$ be an invertible element in $k$
satisfying $q^2\neq 1$. The quantized enveloping algebra\footnote{We
are actually working with De Concini--Kac's version of specialized
quantum algebra, see \cite{CK}.} associated to the simple Lie
algebra $\ksl_2$ is the associative $k$-algebra with the generators
$E, F, K, K^{-1}$ and the relations (cf. \cite{D1} and \cite{J1},
\cite{J2}): $$
\begin{aligned}
& KE=q^2EK,\,\,KF=q^{-2}FK,\,\,KK^{-1}=1=K^{-1}K,\\
& EF-FE=\frac{K-K^{-1}}{q-q^{-1}}.
\end{aligned}
$$
We denote it by ${\ru}_q(\ksl_2)$ or just ${\ru}_q$ for
simplicity. The algebra ${\ru}_q$ is a quantum analogue of the universal
enveloping algebra ${\ru}(\ksl_2)$ associated to the simple Lie algebra $\ksl_2$.
It is a Hopf algebra with
comultiplication, counit and antipode given by:
$$\begin{aligned} &\Delta(E)=E\otimes 1+K\otimes
E,\,\,\Delta(F)=F\otimes K^{-1}+1\otimes
F,\,\,\Delta(K)=K\otimes K, \\
& \varepsilon(E)=0=\varepsilon(F),\,\,\varepsilon(K)=1=\varepsilon(K^{-1}),\\
& S(E)=-K^{-1}E,\,\,S(F)=-FK,\,\,S(K)=K^{-1}.
\end{aligned}
$$
Let ${\ru}_q^{+}$ (resp. ${\ru}_q^{-}$) be the $k$-subalgebra of
${\ru}_q$ generated by $E$ (resp. by $F$). Let ${\ru}_q^{0}$ be
the $k$-subalgebra of ${\ru}_q$ generated by $K, K^{-1}$. Then the
elements $\{E^a\}$ (resp. $\{F^b\}$), where
$a,b\in\mathbb{N}\cup\{0\}$, form a $k$-basis of ${\ru}_q^{+}$
(resp. of ${\ru}_q^{-}$). The elements $\{K^c\}$, where
$c\in\mathbb{Z}$, form a $k$-basis of ${\ru}_q^{0}$. Moreover, the
natural $k$-linear map ${\ru}_q^{+}\otimes{\ru}_q^{0}\otimes
{\ru}_q^{-}\rightarrow {\ru}_q$ given by multiplication is a
$k$-linear isomorphism. The basis
$\bigl\{E^aK^cF^b\bigm|a,b\in\mathbb{N}\cup\{0\},
c\in\mathbb{Z}\bigr\}$ is called the PBW basis of ${\ru}_q$. We
define ${\ru}_q^{\geq 0}:={\ru}_q^{+}{\ru}_q^{0}$, ${\ru}_q^{\leq
0}:={\ru}_q^{-}{\ru}_q^{0}$. Then both ${\ru}_q^{\geq 0}$ and
${\ru}_q^{\leq 0}$ are Hopf $k$-subalgebras of ${\ru}_q$. For any
monomials $E^{a}K^{b}, F^{c}K^{d}$, we endow them {\it the degree
$a, c$} respectively.

Let $v$ be an indeterminate over $\mathbb{Z}$. We consider the
quantized enveloping algebra ${\ru}_v={\ru}_v(\ksl_2)$ with
parameter $v$ and defined over $\mathbb{Q}(v)$. It is well-known
(see \cite{Jo}, \cite{L5} and \cite{X}\footnote{Note that we use a
slightly different version here so that the multiplication rule is
compatible with the one given in \cite[Chapter IX, (4.3)]{K} for
the Drinfel'd quantum double of finite dimensional Hopf
algebras.}) that there exists a unique pairing $\varphi:
{\ru}_v^{\geq 0}\times {\ru}_v^{\leq 0}\rightarrow\mathbb{Q}(v)$
such that
$$\begin{aligned} (1)
&\,\,\,\varphi(1,1)=1,\,\,\varphi(1,K)=1=\varphi(K,1)\\
(2) &\,\,\,\varphi(x,y)=0,\,\,\text{if $x,y$ are homogeneous with different degree,}\\
(3) &\,\,\,\varphi(E,F)=\frac{1}{v^2-1},\,\,\varphi(K,K)=v^{2},\,\,\varphi(K,K^{-1})=v^{-2},\\
(4) &\,\,\,\varphi(x,y'y'')=\varphi(\Delta^{\op}(x),y'\otimes
y''),\,\,\text{for all $x\in{\ru}_v^{\geq 0},\,
y',y''\in{\ru}_v^{\leq 0}$},\\
(5) &\,\,\,\varphi(xx',y'')=\varphi(x\otimes
x',\Delta(y'')),\,\,\text{for all $x,x'\in{\ru}_v^{\geq 0},\,
y''\in{\ru}_v^{\leq 0}$},\\
(6) &\,\,\,\varphi(S(x),y)=\varphi(x,S^{-1}(y)),\,\,\text{for all
$x\in{\ru}_v^{\geq 0},\, y\in{\ru}_v^{\leq 0}$}.
\end{aligned}
$$
One usually call $({\ru}_v^{\geq 0}, {\ru}_v^{\leq 0}, \varphi)$ a
skew Hopf pairing (cf. \cite{M}). Then, one can make
$D({\ru}_v^{\geq 0}, {\ru}_v^{\leq 0}):={\ru}_v^{\geq 0}\otimes
{\ru}_v^{\leq 0}$ into a Hopf $\mathbb{Q}(v)$-algebra, which is
called the quantum double of $({\ru}_v^{\geq 0}, {\ru}_v^{\leq
0},\varphi)$. As a $\mathbb{Q}(v)$-coalgebra, $D({\ru}_v^{\geq 0},
{\ru}_v^{\leq 0})={\ru}_v^{\geq 0}\otimes {\ru}_v^{\leq 0}$, the
tensor product of two coalgebras. The algebra structure of
$D({\ru}_v^{\geq 0}, {\ru}_v^{\leq 0})$ is determined by $$
\begin{aligned}
(7) &\,\,\,(x\otimes y)(x'\otimes
y')=\sum\varphi\bigl(x'_{(1)},y_{(1)}\bigr)(xx'_{(2)}\otimes
y_{(2)}y') \varphi\bigl(x'_{(3)},S^{-1}(y_{(3))}\bigr),
\end{aligned}
$$
for all $x,x'\in {\ru}_v^{\geq 0}, y,y'\in {\ru}_v^{\leq 0}$, where
Sweedler's sigma summation $\Delta^2(y)=\sum y_{(1)}\otimes
y_{(2)}\otimes y_{(3)}$ is used. Note that the quantum double we
described here actually arises from a $2$-cocycle twist (see
\cite{DT} for more detail). The representation results in this paper
are related to the recent work of Radford and Schneider (see
\cite{RS1} and \cite{RS2}). We thank the referee for pointing out
these references. For simplicity, we shall write $\wtu_v$ instead of
$D({\ru}_v^{\geq 0}, {\ru}_v^{\leq 0})$.

Let $A=\mathbb{Z}[v,v^{-1},(v-v^{-1})^{-1}]$. Let ${\ru}_{A,v}^{\leq 0}$ be the associative
$A$-algebra defined by the generators $F, K, K^{-1}$ and relations $$
KF=v^{-2}FK,\,\,KK^{-1}=1=K^{-1}K.
$$
Specializing $v$ to $q$, we make $k$ into an $A$-algebra. Clearly ${\ru}_q^{\leq 0}\cong k\otimes_{A}{\ru}_{A,v}^{\leq 0}$.
Similarly, we can define ${\ru}_{A,v}^{\geq 0}$, and we have ${\ru}_q^{\geq 0}\cong k\otimes_{A}{\ru}_{A,v}^{\geq 0}$.

The previous construction of skew Hopf pairing clearly gives rise to a pairing
$\varphi: {\ru}_{A,v}^{\geq 0}\times {\ru}_{A,v}^{\leq 0}\rightarrow A$, and hence gives rise to a pairing  $\varphi:
{\ru}_q^{\geq 0}\times {\ru}_q^{\leq 0}\rightarrow k$.

\begin{lemma} \lelabel{11} The pairing $\varphi$ gives rise to a Hopf algebra map $\theta$ from the Hopf algebra
${\ru}_q^{\leq 0}$ to the Hopf algebra $\bigl({\ru}_q^{\geq 0}\bigr)^{*,\op}$ as well as a Hopf algebra map $\theta'$
from the Hopf algebra ${\ru}_q^{\geq 0}$ to the Hopf algebra $\bigl({\ru}_q^{\leq 0}\bigr)^{*,\cop}$. Moreover,
$\theta,\theta'$ are injective if
$q$ is not a root of unity.
\end{lemma}

\begin{proof} The maps $\theta, \theta'$ are defined by
$$\begin{aligned}
\theta(y)(x)&=\varphi(x,y),\\
\theta'(x)(y)&=\varphi(x,y),\end{aligned}
$$
for any $x\in {\ru}_{q}^{\geq 0}, y\in {\ru}_{q}^{\leq 0}$.

Now the first statement follows directly from the definition of
$\varphi$. The second statement can be proved by using a similar
argument in the proof of \leref{41}.\end{proof}

Now we can construct the quantum double $D({\ru}_q^{\geq 0},
{\ru}_q^{\leq 0}):={\ru}_q^{\geq 0}\otimes {\ru}_q^{\leq 0}$ of $({\ru}_q^{\geq 0}, {\ru}_q^{\leq 0},\varphi)$ in a
similar way, making it into a Hopf $k$-algebra. Henceforth, we shall write $\wtu_q$ instead of
$D({\ru}_q^{\geq 0},{\ru}_q^{\leq 0})$. We have the following.

\begin{theorem} As a $k$-algebra, $\wtu_q$ can be presented by the generators $$E, F, K, K^{-1},
\widetilde{K}, \widetilde{K}^{-1}, $$ and the following relations:
$$\begin{aligned}
&
KE=q^2EK,\,\,KF=q^{-2}FK,\,\,\widetilde{K}E=q^2E\widetilde{K},\,\,\widetilde{K}F=q^{-2}F\widetilde{K},\\
&KK^{-1}=K^{-1}K=1=\widetilde{K}\widetilde{K}^{-1}=\widetilde{K}^{-1}\widetilde{K},\,\,K\widetilde{K}=\widetilde{K}K,\\
& EF-FE=\frac{K-\widetilde{K}^{-1}}{q-q^{-1}}.
\end{aligned}
$$
\end{theorem}

\begin{proof} Let $\wtu'_q$ be an abstract $k$-algebra defined by the
generators and relations as above. One checks directly that $\wtu_q$
is generated by the following elements $$ E\otimes 1,\,\,1\otimes
qF,\,\,K^{\pm 1}\otimes 1,\,\,1\otimes K^{\pm 1},
$$
and these elements satisfy the above relations. In other words,
there is a surjective $k$-algebra homomorphism $\psi:
\wtu'_q\rightarrow \wtu_q$ such that $$ \psi(E)=E\otimes 1,\,\,
\psi(F)=1\otimes qF,\,\,\psi(K^{\pm 1})=K^{\pm 1}\otimes
1,\,\,\psi(\widetilde{K}^{\pm 1})=1\otimes K^{\pm 1}.$$ On the
other hand, we claim that the monomials $$
E^{a}K^{c}\widetilde{K}^{d}F^{b},\quad
a,b,c,d\in\mathbb{Z},\,a,b\geq 0,
$$
form a basis of $\wtu'_q$. We prove this by using a similar
argument as in the proof of \cite[Theorem 1.5]{Ja}. We consider a
polynomial ring $k[T_1,T_2,T_3,T_4]$ in four indeterminate
$T_1,T_2,T_3,T_4$  and its localization
$A'=k[T_1,T_2,T_2^{-1},T_3,T_{3}^{-1},T_4]$. Then all monomials
$T_1^{a}T_2^{c}T_3^{d}T_4^{b}$ with $a,b,c,d\in\mathbb{Z}, a,b\geq
0$ are a basis of $A'$. We define linear endomorphisms
$e,f,h,\widetilde{h}$ of $A'$  by letting
$$
\begin{aligned}
e\bigl(T_1^{a}T_2^{c}T_3^{d}T_4^{b}\bigr)&=T_1^{a+1}T_2^{c}T_3^{d}T_4^{b},\\
f\bigl(T_1^{a}T_2^{c}T_3^{d}T_4^{b}\bigr)&=\begin{cases}
\begin{matrix}-T_1^{a-1}[a-1]_q\DDF{q^{a-1}K-q^{1-a}\widetilde{K}^{-1}}{q-q^{-1}}T_2^{c}T_3^{d}T_4^{b}\\
+q^{2c+2d}T_1^{a}T_2^{c}T_3^{d}T_4^{b+1},\end{matrix}&\text{if
$a\geq 1$;}\\[2pt]
q^{2c+2d}T_2^{c}T_3^{d}T_4^{b+1}, &\text{if $a=0$,}
\end{cases}
\\
h\bigl(T_1^{a}T_2^{c}T_3^{d}T_4^{b}\bigr)&=q^{2a}T_1^{a}T_2^{c+1}T_3^{d}T_4^{b},\\
\widetilde{h}\bigl(T_1^{a}T_2^{c}T_3^{d}T_4^{b}\bigr)&=q^{2a}T_1^{a}T_2^{c}T_3^{d+1}T_4^{b},
\end{aligned}
$$
where $$[a-1]_q:=\frac{q^{a-1}-q^{1-a}}{q-q^{-1}}.$$ One can check
that the above definition gives rise to a representation of
$\wtu'_q$ on $A'$ by taking $E$ to $e$, $F$ to $f$, $K^{\pm 1}$ to
$h^{\pm 1}$  and $\widetilde{K}^{\pm 1}$ to $\widetilde{h}^{\pm
1}$. So it takes a monomial $E^{a}K^{c}\widetilde{K}^{d}F^{b}$ to
the monomial $e^{a}h^{c}\widetilde{h}^{d}f^{b}$. Note that
$$ e^{a}h^{c}\widetilde{h}^{d}f^{b}(1)=T_1^{a}T_2^{c}T_3^{d}T_4^{b},
$$ which implies that the $e^{a}h^{c}\widetilde{h}^{d}f^{b}$ are
linearly independent, hence the monomials
$E^{a}K^{c}\widetilde{K}^{d}F^{b}$ must  be linear
independent as well. This proves our claim. By definition of $\wtu_q$, we
know that the monomials
$$ (E\otimes 1)^{a}(K\otimes 1) ^{c}(1\otimes
\widetilde{K})^{d}(1\otimes F)^{b},\quad
a,b,c,d\in\mathbb{Z},\,a,b\geq 0,
$$
are a basis of $\wtu_q$. Therefore, $\psi$ maps a basis of
$\wtu'_q$ onto a basis of $\wtu_q$. It follows that $\psi$ must
be an isomorphism, as required.
\end{proof}

\begin{lemma} The map which sends $E$ to $E$, $F$ to $F$, $K^{\pm
1}$ to $K^{\pm 1}$ and $\widetilde{K}^{\pm 1}$ to $K^{\pm 1}$
extends uniquely to a surjective Hopf algebra homomorphism $\pi:
\wtu_q\twoheadrightarrow{\rm U}_q$.
\end{lemma}

\begin{proof} This is obvious. Note that the kernel of $\pi$ is the two-sided
ideal of $\wtu_q$ generated by $K-\widetilde{K}$, which is in fact
a Hopf ideal of $\wtu_q$.
\end{proof}

The algebra $\wtu_q$ will be the primary interest to us in this
paper. It turns out that this algebra behaves quite similar to
the quantized enveloping algebra ${\rm U}_q$ in many ways. In the
following sections we shall see that many constructions and equalities in
the structure and representation theory of ${\rm U}_q$ carry over
to the algebra $\wtu_q$.

\bigskip\bigskip

\section{An analogue of Kostant--Lusztig $\mathbb{Z}[v,v^{-1}]$-form}

The purpose of this section is to construct an analogue of
Kostant--Lusztig $\mathbb{Z}[v,v^{-1}]$-form for the quantum double
${\wtu}_q$.

Let $R$ be an integral domain. Let $v$ be an indeterminate over $R$.
Let $\mathcal{A}=R[v,v^{-1}]$, the ring of Laurent $R$-polynomials
in $v$. We consider the quantum double $\wtu_v$ (with parameter $v$)
defined over the quotient field of $\mathcal{A}$. We shall define an
analogue of Kostant--Lusztig $\mathcal{A}$-form (see \cite{L1},
\cite{L3}, \cite{L4}) for the algebra ${\wtu}_v$. For each positive
integer $N$, we define $$ [N]:=\frac{v^{N}-v^{-N}}{v-v^{-1}},\,\,\,
[N]^{!}:= [N][N-1]\cdots [2][1].
$$
For any integers $m,n$ with $n\geq 0$, we define $$
\begin{bmatrix}m\\
n\end{bmatrix}=\frac{[m]^!}{[n]^![m-n]^!}.
$$ Then it is well-known that $\begin{bmatrix}m\\
n\end{bmatrix}\in\mathcal{A}$ (e.g, see \cite[(1.3.1.d)]{L5}).

Let $\bwtu_{\mathcal{A}}$ be the $\mathcal{A}$-subalgebra of
$\wtu_v$ generated by (compare with \cite[(3.1)]{L1})
$$\begin{aligned} & E^{(N)}=\frac{E^N}{[N]^!},\,\,
F^{(N)}=\frac{F^N}{[N]^!},\,\, K^{\pm 1},\,\,\widetilde{K}^{\pm
1},\\
&\begin{bmatrix}K,\,\widetilde{K}\\
\!\!t\end{bmatrix}=\prod_{s=1}^{t}\frac{Kv^{-s+1}-\widetilde{K}^{-1}v^{s-1}}
{v^{s}-v^{-s}},
\end{aligned}
$$ where $N, t\in\mathbb{N}\cup\{0\}$.

For any integers $c, t$ with $t\geq 0$, we define the analogue of
$\begin{bmatrix}K,c\\
\!\!t\end{bmatrix}$ (see \cite[(4.1)]{L1}).
$$
\begin{bmatrix}K,\widetilde{K},c\\
\!\!t\end{bmatrix}=\prod_{s=1}^{t}\frac{Kv^{c-s+1}-\widetilde{K}^{-1}v^{-c+s-1}}
{v^s-v^{-s}}.$$
We have the following (compare it with \cite[(4.3.1)]{L2}).

\begin{lemma} \lelabel{21} For any non-negative integers $a, b$, we have that
$$
E^{(a)}F^{(b)}=\sum_{0\leq t\leq\min(a,b)}F^{(b-t)}\begin{bmatrix}K,\widetilde{K},2t-a-b\\
\!\!t\end{bmatrix}E^{(a-t)}. $$
\end{lemma}

\begin{proof} The equality is proved in a similar way to the
proof of \cite[(4.3.1)]{L2}. We show it by induction on $a$. For
$a=0$ (or $b=0$), the claim is trivial. For $a=1$ and $b>0$, it
also follows from a straightforward verification. Suppose the
formula holds for integers $a, b$. Then we obtain
$$\begin{aligned}
&\quad\,\,E^{(a+1)}F^{(b)}=\frac{1}{[a+1]}EE^{(a)}F^{(b)}\\
&=\frac{1}{[a+1]}\sum_{0\leq t\leq\min(a,b)}EF^{(b-t)}\begin{bmatrix}K,\widetilde{K},2t-a-b\\
\!\!t\end{bmatrix}E^{(a-t)}\\
&=\frac{1}{[a+1]}\biggl\{\sum_{0\leq t\leq\min(a,b)}F^{(b-t)}E
\begin{bmatrix}K,\widetilde{K},2t-a-b\\
\!\!t\end{bmatrix}E^{(a-t)}+\\
&\quad\sum_{0\leq
t\leq\min(a,b)}F^{(b-t-1)}\frac{Kv^{1-b+t}-\widetilde{K}^{-1}v^{-1+b-t}}{v-v^{-1}}
\begin{bmatrix}K,\widetilde{K},2t-a-b\\
\!\!t\end{bmatrix}E^{(a-t)}\biggr\}\\
&=\frac{1}{[a+1]}\biggl\{\sum_{0\leq
t\leq\min(a,b)}[a-t+1]F^{(b-t)}
\begin{bmatrix}K,\widetilde{K},2t-a-b-2\\
\!\!t\end{bmatrix}E^{(a-t+1)}+\\
&\,\,\sum_{0\leq
t\leq\min(a,b)}F^{(b-t-1)}\frac{Kv^{1-b+t}-\widetilde{K}^{-1}v^{-1+b-t}}{v-v^{-1}}
\begin{bmatrix}K,\widetilde{K},2t-a-b\\
\!\!t\end{bmatrix}E^{(a-t)}\biggr\}.
\end{aligned}
$$
Thus we have to show that the equality
$$
\begin{aligned}
&\quad\,\,\frac{1}{[a+1]}\biggl\{[a-t+1]F^{(b-t)}
\begin{bmatrix}K,\widetilde{K},2t-a-b-2\\
\!\!t\end{bmatrix}E^{(a-t+1)}+\\
&\quad\,\,F^{(b-t)}\frac{Kv^{-b+t}-\widetilde{K}^{-1}v^{b-t}}{v-v^{-1}}
\begin{bmatrix}K,\widetilde{K},2t-2-a-b\\
\!\!t-1\end{bmatrix}E^{(a-t+1)}\biggr\}\\
&=F^{(b-t)}\begin{bmatrix}K,\widetilde{K},2t-a-1-b\\
\!\!t\end{bmatrix}E^{(a+1-t)}
\end{aligned}
$$
holds when $0\leq t\leq\min(a,b)$, and that  the equality
$$\begin{aligned}
&\frac{1}{[a+1]}F^{(b-a-1)}\frac{Kv^{-b+a+1}-\widetilde{K}^{-1}v^{b-a-1}}{v-v^{-1}}
\begin{bmatrix}K,\widetilde{K},a-b\\
\!\!a\end{bmatrix}\\
&=F^{(b-a-1)}\begin{bmatrix}K,\widetilde{K},a+1-b\\
\!\!a+1\end{bmatrix}\end{aligned}
$$
holds when $b\geq a+1$.

The verification of the second equality is straightforward while
the first equality follows from the following equation which can
be calculated easily:
 $$\begin{aligned}
\begin{bmatrix}K,\widetilde{K},2t-a-1-b\\
\!\!t\end{bmatrix}&=\frac{1}{[a+1]}\biggl\{[a-t+1]
\begin{bmatrix}K,\widetilde{K},2t-a-b-2\\
\!\!t\end{bmatrix}+\\
&\quad\,\,\frac{Kv^{-b+t}-\widetilde{K}^{-1}v^{b-t}}{v-v^{-1}}
\begin{bmatrix}K,\widetilde{K},2t-2-a-b\\
\!\!t-1\end{bmatrix}\biggr\}.
\end{aligned}
$$
\end{proof}

Let $\bwtu_{\mathcal{A}}^{+}$ (resp. $\bwtu_{\mathcal{A}}^{-}$) be
the $\mathcal{A}$-subalgebra of $\bwtu_{\mathcal{A}}$ generated by
$E^{(a)}$ (resp. $F^{(b)}$), where $a,b\in\mathbb{N}\cup\{0\}$.
Let $\bwtu_{\mathcal{A}}^{0}$ be the $\mathcal{A}$-subalgebra of
$\bwtu_{\mathcal{A}}$ generated by $
K^{\pm 1},\,\,\widetilde{K}^{\pm 1},\,\,\begin{bmatrix}K,\,\widetilde{K}\\
\!\!t\end{bmatrix},\,\,t=0,1,2,\cdots$. We have the following
analogues of \cite[(2.3),(g8),(g9),(g10))]{L3} and
\cite[(4.1),(d)]{L1})\footnote{Note that in
\cite[(2.3),(g10)]{L3}, the range of $j$ in the summation should
be $0\leq j\leq c$ instead of $0\leq j\leq t$.}

\begin{lemma} \lelabel{22} $\mathrm{1)}$ For any integers $c,t,p$ with $t\geq 0, 0\leq p\leq t$, we have
$$ v^{-pt}\begin{bmatrix}K,\,\widetilde{K},\,c\\
t\end{bmatrix}=\sum_{j=0}^{p}\begin{bmatrix}K,\,\widetilde{K},\,c-p\\
t-j\end{bmatrix}\begin{bmatrix}p\\
j\end{bmatrix}\widetilde{K}^{-j}v^{-cj}.$$ In particular, for any
integer $c,t$ with $0\leq c\leq t$, we have $$
\begin{bmatrix}K,\,\widetilde{K},\,c\\
t\end{bmatrix}=\sum_{j=0}^{c}\begin{bmatrix}K,\,\widetilde{K}\\
t-j\end{bmatrix}\begin{bmatrix}c\\
j\end{bmatrix}\widetilde{K}^{-j}v^{c(t-j)}. $$

$\mathrm{2)}$ For any integers $c,t,p$ with $t\geq 0, p\geq 1$, we have
$$ v^{-pt}\begin{bmatrix}K,\,\widetilde{K},\,-c\\
t\end{bmatrix}=\sum_{j=0}^{t}(-1)^{j}\begin{bmatrix}K,\,\widetilde{K},\,p-c\\
t-j\end{bmatrix}\begin{bmatrix}p+j-1\\
j\end{bmatrix}{K}^{j}v^{-cj}.$$ In particular, for any integer
$c,t$ with $c\geq 1,t\geq 0$, we have $$
\begin{bmatrix}K,\,\widetilde{K},\,-c\\
t\end{bmatrix}=\sum_{j=0}^{t}(-1)^{j}\begin{bmatrix}K,\,\widetilde{K}\\
t-j\end{bmatrix}\begin{bmatrix}c+j-1\\
j\end{bmatrix}{K}^{j}v^{c(t-j)}. $$

$\mathrm{3)}$ For any $c\in\mathbb{Z}, t\in\mathbb{N}\cup\{0\}$,
$\begin{bmatrix}K,\,\widetilde{K},\,c\\
\!\!t\end{bmatrix}\in\bwtu_{\mathcal{A}}^{0}$.

$\mathrm{4)}$ For any non-negative integers $t,t'$ with $t\geq 1$, we have
that
$$
\begin{bmatrix}t+t'\\
t\end{bmatrix}\begin{bmatrix}K,\,\widetilde{K}\\
\!\!t+t'\end{bmatrix}=\sum_{0\leq j\leq
t'}(-1)^{j}v^{t(t'-j)}\begin{bmatrix}t+j-1\\
j\end{bmatrix}K^{j}\begin{bmatrix}K,\,\widetilde{K}\\
\!\!t\end{bmatrix}\begin{bmatrix}K,\,\widetilde{K}\\
\!\!t'-j\end{bmatrix}.$$
\end{lemma}

\begin{proof} 1) The second statement follows from induction on
$t$. It suffices to prove the first statement. We show it by induction on
$p$. The case where $p=0$ is trivial. For $p=1$, we have the following
$$
\begin{aligned}
&\quad\,\,\begin{bmatrix}K,\,\widetilde{K},\,c-1\\
t\end{bmatrix}+\begin{bmatrix}K,\,\widetilde{K},\,c-1\\
t-1\end{bmatrix}\widetilde{K}^{-1}v^{-c}\\
&=\biggl(\prod_{s=1}^{t-1}\frac{Kv^{c-1-s+1}-\widetilde{K}^{-1}v^{-c+1+s-1}}{v^{s}-v^{-s}}\biggr)\biggl\{
\frac{Kv^{c-t}-\widetilde{K}^{-1}v^{-c+t}}{v^{t}-v^{-t}}+\widetilde{K}^{-1}v^{-c}\biggr\}\\
&=v^{-t}\begin{bmatrix}K,\,\widetilde{K},\,c\\
t\end{bmatrix},
\end{aligned}
$$
as required. Suppose now the equality holds for $p=N$, we consider
the case where $p=N+1$. We get
$$\begin{aligned}
&\quad\,\,v^{-(N+1)t}\begin{bmatrix}K,\,\widetilde{K},\,c\\
t\end{bmatrix}=v^{-Nt}v^{-t}\begin{bmatrix}K,\,\widetilde{K},\,c\\
t\end{bmatrix}\\
&=v^{-Nt}\biggl\{\begin{bmatrix}K,\,\widetilde{K},\,c-1\\
t\end{bmatrix}+\begin{bmatrix}K,\,\widetilde{K},\,c-1\\
t-1\end{bmatrix}\widetilde{K}^{-1}v^{-c}\biggr\}\\
&=v^{-Nt}\begin{bmatrix}K,\,\widetilde{K},\,c-1\\
t\end{bmatrix}+v^{-N(t-1)}\begin{bmatrix}K,\,\widetilde{K},\,c-1\\
t-1\end{bmatrix}\widetilde{K}^{-1}v^{-c-N}\\
&=\sum_{j=0}^{N}\biggl\{\begin{bmatrix}K,\,\widetilde{K},\,c-N-1\\
t-j\end{bmatrix}\begin{bmatrix}N\\
j\end{bmatrix}\widetilde{K}^{-j}v^{-(c-1)j}+\\
&\qquad\qquad \begin{bmatrix}K,\,\widetilde{K},\,c-N-1\\
t-1-j\end{bmatrix}\begin{bmatrix}N\\
j\end{bmatrix}\widetilde{K}^{-j-1}v^{-(c-1)j-c-N}\biggr\}\\
&=\sum_{j=0}^{N+1}\begin{bmatrix}K,\,\widetilde{K},\,c-N-1\\
t-j\end{bmatrix}\begin{bmatrix}N+1\\
j\end{bmatrix}\widetilde{K}^{-j}v^{-cj},
\end{aligned}
$$
as desired.

2) Now we use induction on $t$. The case where $t=0$ is trivial. For
$t=1$, we obtain the equations
$$
\begin{aligned}
&\quad\,\,\sum_{j=0}^{1}(-1)^{j}\begin{bmatrix}K,\,\widetilde{K},\,p-c\\
1-j\end{bmatrix}\begin{bmatrix}p+j-1\\
j\end{bmatrix}{K}^{j}v^{-cj}\\
&=\frac{Kv^{p-c}-\widetilde{K}^{-1}v^{c-p}}{v-v^{-1}}-\frac{v^{p}-v^{-p}}{v-v^{-1}}Kv^{-c}\\
&=v^{-p}\frac{Kv^{-c}-\widetilde{K}^{-1}v^{c}}{v-v^{-1}}
=v^{-p}\begin{bmatrix}K,\,\widetilde{K},\,-c\\
1\end{bmatrix},
\end{aligned}
$$
as required.

Suppose the equality holds for $t$, we now consider the equality
for $t+1$. We get then $$\begin{aligned}
&\quad\,\,v^{-p(t+1)}\begin{bmatrix}K,\,\widetilde{K},\,-c\\
t+1\end{bmatrix}=v^{-pt}\begin{bmatrix}K,\,\widetilde{K},\,-c\\
t\end{bmatrix}v^{-p}\frac{Kv^{-c-t}-\widetilde{K}^{-1}v^{c+t}}{v^{t+1}-v^{-t-1}}\\
&=\sum_{j=0}^{t}(-1)^{j}\begin{bmatrix}K,\,\widetilde{K},\,p-c\\
t-j\end{bmatrix}\begin{bmatrix}p+j-1\\
j\end{bmatrix}{K}^{j}v^{-cj-p}\frac{Kv^{-c-t}-\widetilde{K}^{-1}v^{c+t}}{v^{t+1}-v^{-t-1}}\\
&=\begin{bmatrix}K,\,\widetilde{K},\,p-c\\
t+1\end{bmatrix}-\begin{bmatrix}K,\,\widetilde{K},\,p-c\\
t\end{bmatrix}Kv^{-c-t}\frac{v^{p}-v^{-p}}{v^{t+1}-v^{-t-1}}+\\
&\quad\,\,\sum_{j=1}^{t}(-1)^{j}\begin{bmatrix}K,\,\widetilde{K},\,p-c\\
t-j\end{bmatrix}\begin{bmatrix}p+j-1\\
j\end{bmatrix}{K}^{j}v^{-cj-p}\frac{Kv^{-c-t}-\widetilde{K}^{-1}v^{c+t}}{v^{t+1}-v^{-t-1}}.\\
\end{aligned}
$$
Note that $$\begin{aligned}
&\quad
\,\,(-1)^{j}\begin{bmatrix}K,\,\widetilde{K},\,p-c\\
t-j\end{bmatrix}\begin{bmatrix}p+j-1\\
j\end{bmatrix}{K}^{j}v^{-cj-p}\frac{Kv^{-c-t}-\widetilde{K}^{-1}v^{c+t}}{v^{t+1}-v^{-t-1}}\\
&=(-1)^{j}\begin{bmatrix}K,\,\widetilde{K},\,p-c\\
t-j\end{bmatrix}\begin{bmatrix}p+j-1\\
j\end{bmatrix}{K}^{j}v^{-cj}\biggl(\frac{Kv^{-c+p-t+j}-\widetilde{K}^{-1}v^{c-p+t-j}}{v^{t+1-j}-v^{j-t-1}}\\
&\quad\frac{v^{t+1}-v^{-t+2j-1}}{v^{t+1}-v^{-t-1}}-
Kv^{-c-t}\frac{v^{p+2j}-v^{-p}}{v^{t+1}-v^{-t-1}}\biggr)\\
&=(-1)^{j}\begin{bmatrix}K,\,\widetilde{K},\,p-c\\
t+1-j\end{bmatrix}\begin{bmatrix}p+j-1\\
j\end{bmatrix}{K}^{j}v^{-cj}\frac{v^{t+1}-v^{-t+2j-1}}{v^{t+1}-v^{-t-1}}+\\
&\quad\,\,(-1)^{j+1}\begin{bmatrix}K,\,\widetilde{K},\,p-c\\
t-j\end{bmatrix}\begin{bmatrix}p+j-1\\
j\end{bmatrix}{K}^{j+1}v^{-c(j+1)}
v^{-t}\frac{v^{p+2j}-v^{-p}}{v^{t+1}-v^{-t-1}}.
\end{aligned}
$$
Now the required equality follows from the following calculation:
$$\begin{aligned} &\quad\,\,(-1)^{j}\begin{bmatrix}K,\,\widetilde{K},\,p-c\\
t+1-j\end{bmatrix}\begin{bmatrix}p+j-1\\
j\end{bmatrix}{K}^{j}v^{-cj}\frac{v^{t+1}-v^{-t+2j-1}}{v^{t+1}-v^{-t-1}}+\\
&\quad\,\,(-1)^{j}\begin{bmatrix}K,\,\widetilde{K},\,p-c\\
t-j+1\end{bmatrix}\begin{bmatrix}p+j-2\\
j-1\end{bmatrix}{K}^{j}v^{-cj}
v^{-t}\frac{v^{p+2j-2}-v^{-p}}{v^{t+1}-v^{-t-1}}\\
&=(-1)^{j}\begin{bmatrix}K,\,\widetilde{K},\,p-c\\
t+1-j\end{bmatrix}\begin{bmatrix}p+j-1\\
j\end{bmatrix}{K}^{j}v^{-cj}\biggl(\frac{v^{t+1}-v^{-t+2j-1}}{v^{t+1}-v^{-t-1}}+\\
&\quad\,\,\frac{[j]}{[p+j-1]}\frac{v^{p+2j-t-2}-v^{-p-t}}{v^{t+1}-v^{-t-1}}\biggr)\\
&=(-1)^{j}\begin{bmatrix}K,\,\widetilde{K},\,p-c\\
t+1-j\end{bmatrix}\begin{bmatrix}p+j-1\\
j\end{bmatrix}{K}^{j}v^{-cj}.\\
\end{aligned}
$$

3) follows from 1) and 2).

4) follows from 2) and the following equality: $$
\begin{bmatrix}K,\,\widetilde{K}\\
t+t'\end{bmatrix}=\begin{bmatrix}K,\,\widetilde{K}\\
t\end{bmatrix}\begin{bmatrix}K,\,\widetilde{K},\,-t\\
t'\end{bmatrix}.
$$
\end{proof}

Let $\theta$ be the algebra automorphism of
$\wtu_{v}$ which is defined on generators by $$
\theta(E)=E,\quad\theta(F)=F,\quad \theta(K^{\pm
1})=\widetilde{K}^{\pm 1},\quad \theta(\widetilde{K}^{\pm
1})=K^{\pm 1}.
$$
Since $$ \theta\Bigl(\begin{bmatrix}K,\,\widetilde{K}\\
\!\!t\end{bmatrix}\Bigr)=\prod_{s=1}^{t}\frac{\widetilde{K}v^{-s+1}-{K}^{-1}v^{s-1}}{v^{s}-v^{-s}}=
K^{-t}\widetilde{K}^{t}\begin{bmatrix}K,\,\widetilde{K}\\
\!\!t\end{bmatrix},
$$
it follows that $\theta$ restricts to an $\mathcal{A}$-algebra
automorphism of $\bwtu_{\mathcal{A}}$. Henceforth, we write $$
\begin{bmatrix}\widetilde{K},\,K\\
\!\!t\end{bmatrix}:=\theta\Bigl(\begin{bmatrix}K,\,\widetilde{K}\\
\!\!t\end{bmatrix}\Bigr),\quad
\begin{bmatrix}\widetilde{K},\,K,\,c\\
\!\!t\end{bmatrix}:=\theta\Bigl(\begin{bmatrix}K,\,\widetilde{K},\,c\\
\!\!t\end{bmatrix}\Bigr).
$$
Then one can get a second version of our previous two lemmas by
applying the automorphism $\theta$.

\begin{lemma} With the notations as above, we have that
the $\mathcal{A}$-algebra $\bwtu_{\mathcal{A}}^{+}$ (resp.
$\bwtu_{\mathcal{A}}^{-}$) is a free $\mathcal{A}$-module, and the
set $\bigl\{E^{(a)}\bigr\}$ (resp. the set
$\bigl\{F^{(b)}\bigr\}$), where $a, b\in\mathbb{N}\cup\{0\}$, form
an $\mathcal{A}$-basis of $\bwtu_{\mathcal{A}}^{+}$ (resp. of
$\bwtu_{\mathcal{A}}^{-}$),
\end{lemma}

\begin{proof} This follows from the fact that the
set $\bigl\{E^{(a)}\bigr\}$ (resp. the set
$\bigl\{F^{(a)}\bigr\}$) is a basis of $\bwtu_v^{+}$ (resp.
of $\bwtu_v^{-}$) and the following equalities: $$ E^{(a)}E^{(b)}=\begin{bmatrix}a+b\\
\!b\end{bmatrix}E^{(a+b)},\quad F^{(a)}F^{(b)}=\begin{bmatrix}a+b\\
\!b\end{bmatrix}F^{(a+b)}. $$
\end{proof}

We define $\bwtu_{\mathcal{A}}^{\geq
0}:=\bwtu_{\mathcal{A}}^{+}\bwtu_{\mathcal{A}}^{0}$,
$\bwtu_{\mathcal{A}}^{\leq
0}:=\bwtu_{\mathcal{A}}^{-}\bwtu_{\mathcal{A}}^{0}$. According to
\cite[Remark 3.1]{Ja}, we know that for each positive integer $N$,
$$\begin{aligned}
\Delta(E^{(N)})&=\sum_{i=0}^{N}v^{i(N-i)}E^{(N-i)}K^i\otimes
E^{(i)},\\
\Delta(F^{(N)})&=\sum_{i=0}^{N}v^{i(N-i)}F^{(i)}\otimes
F^{(N-i)}\widetilde{K}^{-i},\\
S(E^{(N)})&=(-1)^{N}v^{(1-N)N}K^{-N}E^{(N)},\\
S(F^{(N)})&=(-1)^{N}v^{(N-1)N}F^{(N)}\widetilde{K}^{N}.
\end{aligned}
$$

\begin{lemma} For any positive integer $t$, we have $$
\Delta\Bigl(\begin{bmatrix}K,\,\widetilde{K}\\
\!\!t\end{bmatrix}\Bigr)=\sum_{a=0}^{t}\begin{bmatrix}K,\,\widetilde{K}\\
\!t-a\end{bmatrix}\widetilde{K}^{-a}\otimes K^{t-a}\begin{bmatrix}K,\,\widetilde{K}\\
\!\!a\end{bmatrix}.
$$
\end{lemma}

\begin{proof} We use induction on $t$. If $t=1$, we have that
$$\begin{aligned}
&\quad\,\Delta\Bigl(\begin{bmatrix}K,\,\widetilde{K}\\
\!\!1\end{bmatrix}\Bigr)\\
&=\Delta\Bigl(\frac{K-\widetilde{K}^{-1}}{v-v^{-1}}\Bigr)\\
&=\frac{K\otimes K-\widetilde{K}^{-1}\otimes\widetilde{K}^{-1}}{v-v^{-1}}\\
&=\frac{K-\widetilde{K}^{-1}}{v-v^{-1}}\otimes
K+\widetilde{K}^{-1}\otimes
\frac{K-\widetilde{K}^{-1}}{v-v^{-1}}\\
&=\begin{bmatrix}K,\,\widetilde{K}\\
\!t-a\end{bmatrix}\otimes K+\widetilde{K}^{-1}\otimes\begin{bmatrix}K,\,\widetilde{K}\\
\!t-a\end{bmatrix},
\end{aligned}$$
as required.

Now assume that the equality holds for $t=N$. We consider the case
where $t=N+1$. We have that $$\begin{aligned}
&\quad\,\,\Delta\Bigl(\begin{bmatrix}K,\,\widetilde{K}\\
\!N+1\end{bmatrix}\Bigr)\\
&=\Delta\Bigl(\prod_{s=1}^{N+1}\frac{Kv^{-s+1}-\widetilde{K}^{-1}v^{s-1}}
{v^{s}-v^{-s}}\Bigr)\\
&=\Delta\Bigl(\begin{bmatrix}K,\,\widetilde{K}\\
\!N\end{bmatrix}\Bigr)\Delta\Bigl(\frac{Kv^{-N}-\widetilde{K}^{-1}v^{N}}
{v^{N+1}-v^{-N-1}}\Bigr)\\
&=\sum_{a=0}^{N}\Bigl(\begin{bmatrix}K,\,\widetilde{K}\\
\!N-a\end{bmatrix}\widetilde{K}^{-a}\otimes K^{N-a}\begin{bmatrix}K,\,\widetilde{K}\\
\!\!a\end{bmatrix}\Bigr)\Bigl(\frac{Kv^{-N}\otimes K
-\widetilde{K}^{-1}v^{N}\otimes \widetilde{K}^{-1}}
{v^{N+1}-v^{-N-1}}\Bigr)\\
&=\sum_{a=0}^{N}\Bigl(\begin{bmatrix}K,\,\widetilde{K}\\
\!N-a\end{bmatrix}\widetilde{K}^{-a}\otimes K^{N-a}\begin{bmatrix}K,\,\widetilde{K}\\
\!\!a\end{bmatrix}\Bigr)\Bigl(\frac{Kv^{-N}-\widetilde{K}^{-1}v^{N}}{v^{N+1}-v^{-N-1}}\otimes K+\\
&\qquad\qquad \frac{v^{N}\widetilde{K}^{-1}}{[N+1]}\otimes\begin{bmatrix}K,\,\widetilde{K}\\
\!1\end{bmatrix}\Bigr)\\
&=\sum_{a=0}^{N}\biggl(\begin{bmatrix}K,\,\widetilde{K}\\
\!N-a\end{bmatrix}\widetilde{K}^{-a}\frac{Kv^{-N}-\widetilde{K}^{-1}v^{N}}{v^{N+1}-v^{-N-1}}
\otimes K^{N+1-a}\begin{bmatrix}K,\,\widetilde{K}\\
\!\!a\end{bmatrix}+\\
&\qquad\qquad \frac{v^{N}}{[N+1]}\begin{bmatrix}K,\,\widetilde{K}\\
\!N-a\end{bmatrix}\widetilde{K}^{-a-1}\otimes K^{N-a}\begin{bmatrix}K,\,\widetilde{K}\\
\!\!a\end{bmatrix}\begin{bmatrix}K,\,\widetilde{K}\\
\!1\end{bmatrix}\biggr)\end{aligned} $$

$$\begin{aligned}
&=\sum_{a=0}^{N}\Biggl(\begin{bmatrix}K,\,\widetilde{K}\\
\!N-a\end{bmatrix}\widetilde{K}^{-a}\frac{Kv^{-N}-\widetilde{K}^{-1}v^{N}}{v^{N+1}-v^{-N-1}}
\otimes K^{N+1-a}\begin{bmatrix}K,\,\widetilde{K}\\
\!\!a\end{bmatrix}+\\
&\qquad\qquad \frac{v^{N}}{[N+1]}\begin{bmatrix}K,\,\widetilde{K}\\
\!N-a\end{bmatrix}\widetilde{K}^{-a-1}\otimes \biggl(v^{-a}[a]K^{N+1-a}\begin{bmatrix}K,\,\widetilde{K}\\
\!\!a\end{bmatrix}+\\
&\qquad\qquad\qquad\qquad v^{-a}[a+1]K^{N-a}\begin{bmatrix}K,\,\widetilde{K}\\
\!\!a+1\end{bmatrix}\biggr)\Biggr)\qquad {\rm (\text{by
\leref{22}})}\\
&=\widetilde{K}^{-N-1}\otimes\begin{bmatrix}K,\,\widetilde{K}\\
\!\!N+1\end{bmatrix}+\sum_{a=0}^{N}\Biggl\{\biggl(\begin{bmatrix}K,\,\widetilde{K}\\
\!N-a\end{bmatrix}\widetilde{K}^{-a}\frac{Kv^{-N}-\widetilde{K}^{-1}v^{N}}{v^{N+1}-v^{-N-1}}
+\\
&\qquad \begin{bmatrix}K,\,\widetilde{K}\\
\!N-a\end{bmatrix}\frac{v^{N-a}[a]}{[N+1]}\widetilde{K}^{-a-1}+
\frac{v^{N-a+1}(v^{a}-v^{-a})}{v^{N+1}-v^{-N-1}}
\begin{bmatrix}K,\,\widetilde{K}\\
\!N-a+1\end{bmatrix}\biggr)\otimes\\
&\qquad\qquad K^{N+1-a}\begin{bmatrix}K,\,\widetilde{K}\\
\!\!a\end{bmatrix}\Biggr\}\\
&=\sum_{a=0}^{N+1}\begin{bmatrix}K,\,\widetilde{K}\\
\!N+1-a\end{bmatrix}\widetilde{K}^{-a}\otimes K^{N+1-a}\begin{bmatrix}K,\,\widetilde{K}\\
\!\!a\end{bmatrix}.
\end{aligned}
$$
This proves the lemma.
\end{proof}

Note that a special case of the above result appeared in the proof
of \cite[Lemma 1.1(ii)]{APW}. However, we did not find any
specific reference to the above calculations.

\begin{corollary} With the
notations as above, $\bwtu_{\mathcal{A}}$ is a Hopf
$\mathcal{A}$-subalgebra of $\bwtu_{v}$, and both
$\bwtu_{\mathcal{A}}^{\geq 0}$ and $\bwtu_{\mathcal{A}}^{\leq 0}$
are Hopf subalgebras of $\bwtu_{\mathcal{A}}$.
\end{corollary}

It is easy to see that $\bwtu_{\mathcal{A}}^{\geq
0}\cong\bwtu_{\mathcal{A}}^{+}\otimes\bwtu_{\mathcal{A}}^{0}$,
$\bwtu_{\mathcal{A}}^{\leq
0}\cong\bwtu_{\mathcal{A}}^{-}\otimes\bwtu_{\mathcal{A}}^{0}$. For
any field $k$ which is an $\mathcal{A}$-algebra, we define
$\bwtu_k:=k\otimes_{\mathcal{A}}\bwtu_{\mathcal{A}}$,
$\bwtu_k^{\geq 0}:=k\otimes_{\mathcal{A}}\bwtu_{\mathcal{A}}^{\geq
0}$, $\bwtu_k^{\leq
0}:=k\otimes_{\mathcal{A}}\bwtu_{\mathcal{A}}^{\leq 0}$.

\begin{remark} It would be interesting to know if $\bwtu_{\mathcal{A}}^{0}$ is a free $\mathcal{A}$-module
and whether there is a triangular decomposition for the $\mathcal{A}$-algebra $\bwtu_{\mathcal{A}}$.
\end{remark}

\begin{corollary} We consider $\mathbb{Q}(v)$ as an $\mathcal{A}$-algebra in a natural way. Then the natural map
$\mathbb{Q}(v)\otimes_{\mathcal{A}}\bwtu_{\mathcal{A}}\rightarrow\wtu_{v},\,\,
a\otimes x\mapsto ax,\,\forall\,a\in\mathbb{Q}(v), x\in \bwtu_{\mathcal{A}}$,
is a $\mathbb{Q}(v)$-algebra isomorphism.
\end{corollary}

\bigskip\bigskip

\section{Representations of the algebra $\wtu_q$}

Throughout this section, we assume that $k$ is an algebraically closed field.

Let $M$ be a $\wtu_q$-module such that $\End_{\wtu_q}(M)=k$. Note
that the element $K\widetilde{K}^{-1}$ is invertible and central in
$\wtu_q$. Therefore, there is an element $0\neq z\in k$ such that
$K\widetilde{K}^{-1}$ acts as the scalar $z$ on $M$. For each $0\neq
z\in k$, we fix a square root $\sqrt{z}$ of $z$. Let $\pi_z^{+}$ be
the $k$-algebra homomorphism $\wtu_q\rightarrow{\ru}_q$ which is
defined on generators by $$
\pi_z^{+}(E)=\sqrt{z}E,\,\,\pi_z^{+}(F)=F,\,\,\pi_z^{+}(K)=\sqrt{z}K,\,\,\pi_z^{+}(\widetilde{K})=\sqrt{z}^{-1}K.
$$
It is easy to check that $\pi_z^{+}$ is well-defined. Moreover, the
kernel of $\pi_z^{+}$, which is the ideal generated by
$K\widetilde{K}^{-1}-z$, annihilates the module $M$. It follows that
$M$ becomes a module over the algebra ${\ru}_q$ in a natural way.
Similarly, we have a well-defined $k$-algebra homomorphism
$\pi_z^{-}:\wtu_q\rightarrow{\ru}_q$ which is defined on generators
by $$
\pi_z^{-}(E)=-\sqrt{z}E,\,\,\pi_z^{-}(F)=F,\,\,\pi_z^{-}(K)=-\sqrt{z}K,\,\,\pi_z^{-}(\widetilde{K})=-\sqrt{z}^{-1}K.
$$
Note that $\pi_1^{+}$ is a Hopf algebra map, but in general, both
$\pi_z^{+}$ and $\pi_z^{-}$ are not Hopf algebra maps.

We call a $\wtu_q$-module $M$ a weight $\wtu_q$-module if both $K$
and $\widetilde{K}$ act semisimply on $M$. In that case,
$K\widetilde{K}^{-1}$ acts semisimply on $M$ as well. Similarly, we
call a ${\ru}_q$-module $N$ a weight ${\ru}_q$-module if $K$ acts
semisimply on $N$.

\begin{lemma}\lelabel{31} Every finite dimensional simple (resp. indecomposable weight)
$\wtu_q$-module is the pull-back of a finite dimensional simple
(resp. indecomposable weight) ${\ru}_q$-module through the algebra
homomorphisms $\pi_z^{\pm}$ for some $0\neq z\in k$.
\end{lemma}

\begin{proof} For any finite dimensional $\wtu_q$-module $M$, we consider the
 eigenspace of $K\widetilde{K}^{-1}$ on $M$. Since $K\widetilde{K}^{-1}$ is central in $\wtu_q$,
we deduce that each such eigenspace must be a $\wtu_q$-submodule of
$M$. Therefore, if $M$ is a simple $\wtu_q$-module or an
indecomposable weight $\wtu_q$-module, then $K\widetilde{K}^{-1}$
can have only one eigenvalue on $M$. This proves that
$K\widetilde{K}^{-1}$ acts as a scalar on $M$, hence the lemma
follows immediately from the previous discussion.
\end{proof}

Let $M$ be a ${\ru}_q$-module. For any $0\neq \lam\in k$. We denote
by $M_{\lam}^{+}$ (resp. $M_{\lam}^{-}$) the pull-back of $M$
through the algebra homomorphism $\pi_{\lam}^{+}$ (resp.
$\pi_{\lam}^{-}$).

\begin{theorem} \thlabel{32} The category $\widetilde{\mathcal{C}}$ of finite dimensional weight $\wtu_q$-modules is
equivalent to a direct sum of $|k^{\times}|$ copies of the category
$\mathcal{C}$  of finite dimensional weight ${\ru}_q$-modules.
\end{theorem}

\begin{proof} By definition, every object $M$ of $\mathcal{C}$ is of the
form $\oplus_{\lam\in k^{\times}}M(\lam)$ , where for each $\lam$,
$M(\lam)$ is a finite dimensional indecomposable weight
${\ru}_q$-module, and $|\bigl\{\lam\in k^{\times}\bigm|M(\lam)\neq
0\bigr\}|<\infty$. We use $\theta^{+}$ to denote the functor from
$\mathcal{C}$ to $\widetilde{\mathcal{C}}$ such that $$
\theta^{+}\Bigl(\bigoplus_{\lam\in k^{\times}}M(\lam)\Bigr):=
\bigoplus_{\lam\in k^{\times}}M(\lam)_{\lam}^{+}.
$$
The action of $\theta^{+}$ on the set of morphisms is defined in an
obvious way. Then applying \leref{31}, we see that $\theta^{+}$ is
an equivalence of categories. In a similar way, if we define
$\theta^{-}$ to be the functor from $\mathcal{C}$ to
$\widetilde{\mathcal{C}}$ satisfying $$
\theta^{-}\Bigl(\bigoplus_{\lam\in k^{\times}}M(\lam)\Bigr):=
\bigoplus_{\lam\in k^{\times}}M(\lam)_{\lam}^{-},
$$
and the action of $\theta^{-}$ on the set of morphisms is defined in
an obvious way, then $\theta^{-}$ is also an equivalence of
categories.
\end{proof}

Let $M$ be a ${\ru}_q$-module. Let $0\neq z\in k$. Let
$\varepsilon_z^{+}$ (resp. $\varepsilon_z^{-}$) be the
one-dimensional representation of $\wtu_q$ which is defined on
generators by $$
\varepsilon_z^{+}(E)=0=\varepsilon_z^{+}(F),\,\,\varepsilon_z^{+}(K)=\sqrt{z},\,\,
\varepsilon_z^{+}(\widetilde{K})=\sqrt{z}^{-1}.
$$
$\biggl($ resp. $
\varepsilon_z^{-}(E)=0=\varepsilon_z^{-}(F),\,\,\varepsilon_z^{-}(K)=-\sqrt{z},\,\,
\varepsilon_z^{-}(\widetilde{K})=-\sqrt{z}^{-1}. \biggr)$ It is easy
to check that both $\varepsilon_z^{+}$ and $\varepsilon_z^{-}$ are
well-defined. For any $z,z'\in k^{\times}$, we have that
$$\begin{aligned}
&\varepsilon_z^{\pm}\otimes\varepsilon_{z'}^{\pm}\cong
\varepsilon_{z'}^{\pm}\otimes\varepsilon_{z}^{\pm}\cong\begin{cases}
\varepsilon_{zz'}^{+}, &\text{if $\sqrt{z}\sqrt{z'}=\sqrt{zz'}$;}\\
\varepsilon_{zz'}^{-}, &\text{if $\sqrt{z}\sqrt{z'}=-\sqrt{zz'}$,}
\end{cases}
\\
&\varepsilon_z^{+}\otimes\varepsilon_{z'}^{-}\cong
\varepsilon_{z'}^{-}\otimes\varepsilon_{z}^{+} \cong\begin{cases}
\varepsilon_{zz'}^{+}, &\text{if $\sqrt{z}\sqrt{z'}=-\sqrt{zz'}$;}\\
\varepsilon_{zz'}^{-}, &\text{if $\sqrt{z}\sqrt{z'}=\sqrt{zz'}$.}
\end{cases}
\end{aligned}
$$

\begin{lemma} \lelabel{tensor} Let $0\neq z\in k$, let $M$ be a ${\ru}_q$-module. Then

1) $M_{z}^{+}\cong \varepsilon_z^{+}\otimes M_{1}$;

2)  $\varepsilon_z^{+}\otimes M_{1}\cong
M_1\otimes\varepsilon_z^{+}$ if and only if
$\varepsilon_z^{-}\otimes M_{1}\cong M_1\otimes\varepsilon_z^{-}$,
in that case, for any $0\neq z'\in k$ and any ${\ru}_q$-module $N$,
we have that
$$\begin{aligned} M_z^{\pm}\otimes N_{z'}^{\pm}&\cong\begin{cases}
 (M\otimes N)_{zz'}^{+}, &\text{if
 $\sqrt{z}\sqrt{z'}=\sqrt{zz'}$;}\\
(M\otimes N)_{zz'}^{-}, &\text{if
 $\sqrt{z}\sqrt{z'}=-\sqrt{zz'}$;}\end{cases}\\
M_z^{\pm}\otimes N_{z'}^{\mp}&\cong\begin{cases}
 (M\otimes N)_{zz'}^{+}, &\text{if
 $\sqrt{z}\sqrt{z'}=-\sqrt{zz'}$;}\\
(M\otimes N)_{zz'}^{-}, &\text{if
 $\sqrt{z}\sqrt{z'}=\sqrt{zz'}$;}\end{cases}
\end{aligned} $$

3) $\varepsilon_z^{\pm}\otimes M_{1}\cong
M_1\otimes\varepsilon_z^{\pm}$ if and only if
$\varepsilon_z^{\pm}\otimes \theta^{+}(M)\cong
\theta^{+}(M)\otimes\varepsilon_z^{\pm}$. The same is true if we
replace ``$\theta^{+}$" by ``$\theta^{-}$".

\end{lemma}

\begin{proof} The first statement follows from direct verification. The second and the third statements follow from the
associativity of the tensor product and the previous discussion.
\end{proof}

Note that in general, the assumption $\varepsilon_z^{+}\otimes
M_{1}\cong M_1\otimes\varepsilon_z^{+}$ in \leref{tensor} (2) may
not hold. However, it does hold in the following two
cases:\smallskip

\noindent {\it Case 1.} $q$ is not a root of unity, $M$ is an
integrable weight module over $\ru_q$, i.e., both $E, F$ act locally
nilpotently on $M$ and both $K, {K}^{-1}$ act semisimply on $M$. In
this case, we claim that the $\wtu_q$-modules
$\varepsilon_z^{+}\otimes M_{1}$ and $M_1\otimes\varepsilon_z^{+}$
are isomorphic to each other. In fact, since every integrable weight
module over $\ru_q$ is completely reducible, we can reduce the proof
to the case where $M$ is an irreducible highest weight module over
$\ru_q$. Then $\varepsilon_z^{+}\otimes M_{1}\cong M_z^{+}$ is an
irreducible $\wtu_q$-module. Note that the central element
$K\widetilde{K}^{-1}$ acts as the same scalar on both
$\varepsilon_z^{+}\otimes M_{1}$ and $M_1\otimes\varepsilon_z^{+}$.
It follows that the $\wtu_q$-action on both of these two modules can
factored through the surjective homomorphism $\pi_{z'}^{+}$ for some
$z'\in k^{\times}$. So these two modules can be naturally regarded
as integrable weight modules over $\ru_q$. On the other hand, it is
well-known that the isomorphism class of an integrable weight module
over $\ru_q$ is completely determined by its characters. Since both
$M_1\otimes\varepsilon_z^{+}$ and $\varepsilon_z^{+}\otimes M_1$
have the same characters, they must be isomorphic to each other as
$\ru_q$-modules, and hence are also isomorphic to each other as
$\wtu_q$-modules.

\noindent {\it Case 2.} $q$ is a primitive $d$th root of unity,
$(\sqrt{z})^{d}=1$, $M$ is a weight $\wtu_q$-module such that all of
the elements $E^{d},\,\,F^{d},\,\,K^{d}-1$ act as $0$ on $M$. In
this case, both $\varepsilon_z^{+}$ and $M_1$ can be regarded as
modules over the quotient algebra $$ D_q/\langle E^d, F^d, K^d-1,
\widetilde{K}^d-1\rangle.
$$
Note that the algebra $D_q/\langle E^d, F^d, K^d-1,
\widetilde{K}^d-1\rangle$ is actually a Hopf algebra, and the
natural homomorphism from $D_q$ onto $D_q/\langle E^d, F^d, K^d-1,
\widetilde{K}^d-1\rangle$ is indeed a Hopf algebra homomorphism, and
$D_q/\langle E^d, F^d, K^d-1, \widetilde{K}^d-1\rangle$ is indeed
isomorphic to the Drinfel'd double of a Taft algebra. Hence it is a
quasi-triangular Hopf algebra. As a consequence, both
$M_1\otimes\varepsilon_z^{+}$ and $\varepsilon_z^{+}\otimes M_1$ are
isomorphic to each other as modules over $D_q/\langle E^d, F^d,
K^d-1, \widetilde{K}^d-1\rangle$, and hence are also isomorphic to
each other as $\wtu_q$-modules.
\smallskip

We use $\widetilde{\mathcal{C}}_0$ to denote the full subcategory of
all the finite dimensional weight $\wtu_q$-modules $\widetilde{M}$
satisfying $\varepsilon_z^{+}\otimes \widetilde{M}\cong
\widetilde{M}\otimes\varepsilon_z^{+}$ for any $z\in k^{\times}$,
and we use $\mathcal{C}_0$ to denote the full subcategory of all the
finite dimensional weight $\ru_q$-modules ${M}$ satisfying
$\varepsilon_z^{+}\otimes {M}_1\cong {M}_1\otimes\varepsilon_z^{+}$
for any $z\in k^{\times}$.

\leref{tensor} (2) provides a very easy solution to the problem of
decomposing the tensor product of certain $\wtu_q$-modules, i.e.,
reducing them to the corresponding problem for ${\ru}_q$-modules,
where it has been extensively studied and the results are
well-known, see \cite{P}, \cite{S} and \cite{X2}. Therefore, a large
part of the representations (including all irreducible
representations) of the quantum double $\wtu_q$ can be realized as
certain pullback from the representations of the quantized
enveloping algebra ${\ru}_q$. Note that the representations of the
quantized enveloping algebra ${\ru}_q$ is well-understood (cf.
\cite{Ja}). In particular, the tensor product of finite dimensional
simple $\wtu_q$-modules is determined.

In the following, we shall summarize some results and corollaries
for the algebra $\wtu_q$. We mainly follow the formulation given in
\cite{Ja}. We fix a $0\neq z\in k$. For each $0\neq\lam\in k$, let
$M(\lam)$ be the ${\ru}_q$-module defined in \cite[(2.4)]{Ja}. By
pulling back through $\pi_z^{\pm}$, we get a $\wtu_q$-module
$M_z^{\pm}(\lam)$. We call it the Verma modules associated to
$(z,\lam)$. We have the following two results concerning Verma
modules and simple modules over $D_q$ (compare with \cite[(2.4),
(2.5)]{Ja}).

\begin{corollary}\colabel{34} With the notations as above, $$
M_z^{+}(\lam)\cong\wtu_q/(\wtu_qE+\wtu_q(K-\sqrt{z}\lam)+\wtu_q(\widetilde{K}-\sqrt{z}^{-1}\lam)),
$$
and there is a $k$-basis $\{m_i\}_{i=o}^{\infty}$ of $M_{z}(\lam)$ such that for all $i$,
$$\begin{aligned}
Km_i&=\sqrt{z}\lam q^{-2i}m_i,\,\,\widetilde{K}m_i=\sqrt{z}^{-1}\lam q^{-2i}m_i,\\
Fm_i&=m_{i+1},\\
Em_i&=\begin{cases} 0, &\text{if $i=0$,}\\
[i]_q\sqrt{z}\DDF{\lam q^{1-i}-\lam^{-1}q^{i-1}}{q-q^{-1}}m_{i-1},
&\text{otherwise,}
\end{cases}
\end{aligned}
$$
where $$[i]_q:=\frac{q^i-q^{-i}}{q-q^{-1}}.$$ The result for
$M_z^{-}(\lam)$ is similar.
\end{corollary}

\begin{corollary}\colabel{35} Suppose that $q$ is not a root of unity in $k$ and $0\neq \lam\in k$.
If $\lam\neq \pm q^{n}$ for all integers $n\geq 0$, then the
$\wtu_q$-module $M_z^{\pm}(\lam)$ is simple. If $\lam=\pm q^{n}$ for
some integers $n\geq 0$, then $M_z^{\pm}(\lam)$ has a unique maximal
proper submodule which is spanned by all $m_i$ with $i\geq n+1$ and
is isomorphic to $M_z^{\pm}(q^{-2(n+1)}\lam)$. In this case, the
quotient of $M_z^{\pm}(\lam)$ modulo the maximal proper submodule is
an $(n+1)$-dimensional simple $\wtu_q$-module.
\end{corollary}

Suppose that $q$ is not a root of unity in $k$. By Corollary 3.5, we
know that if $\lam=q^{n}$ for some integers $n\geq 0$, we get two
$(n+1)$-dimensional simple $\wtu_q$-modules, we denote it by
$L_z^{+}(n,+), L_z^{-}(n,+)$; while if $\lam=-q^{n}$ for some
integers $n\geq 0$, we get two $(n+1)$-dimensional simple
$\wtu_q$-modules, denoted by $L_z^{+}(n,-), L_z^{-}(n,-)$. Note that
$L_z^{+}(n,+)\cong L_z^{-}(n,-)$, $L_z^{-}(n,+)\cong L_z^{+}(n,-)$.
Therefore, we define $$
L_z(n,+):=L_z^{+}(n,+),\,\,\,L_z(n,-):=L_z^{+}(n,-).
$$
In fact, the simple $\wtu_q$-module $L_z(n,+)$ (resp. $L_z(n,-)$) is
just the pull-back of simple ${\ru}_q$-module $L(n,+)$ (resp.
$L(n,-)$) through the $k$-algebra homomorphism $\pi_z^{+}$, see
\cite[Theorem 2.6]{Ja} for the definitions of $L(n,+)$ and $L(n,-)$.

By construction, $L_z(n,+)$ has a basis
$\{m_i\}_{i=0}^{n}$ such that $$
\begin{aligned}
Km_i&=z^{1/2}q^{n-2i}m_i,\,\,\widetilde{K}m_i=z^{-1/2}q^{n-2i}m_i,\\
Fm_i&=m_{i+1},\\
Em_i&=\begin{cases} 0, &\text{if $i=0$,}\\
z^{1/2}[i]_q[n+1-i]_q m_{i-1}, &\text{otherwise.}
\end{cases}
\end{aligned}
$$
Similarly, $L_z(n,-)$ has a basis
$\{m'_i\}_{i=0}^{n}$ such that $$
\begin{aligned}
Km'_i&=-z^{1/2}q^{n-2i}m'_i,\,\,\widetilde{K}m'_i=-z^{-1/2}q^{n-2i}m'_i,\\
Fm'_i&=m'_{i+1},\\
Em'_i&=\begin{cases} 0, &\text{if $i=0$,}\\
-z^{1/2}[i]_q[n+1-i]_q m'_{i-1}, &\text{otherwise.}
\end{cases}
\end{aligned}
$$
Note that $L_z(n,+)\not\cong L_z(n,-)$. In fact, $L_z(n,+)\cong
\varepsilon_1^{-}\otimes L_z(n,-)$, where $\varepsilon_1^{-}$ is the
one-dimensional representation of $\wtu_q$ which is defined on
generators by $
\varepsilon_1^{-}(E)=0=\varepsilon_1^{-}(F),\,\,\varepsilon_1^{-}(K)=-1=\varepsilon_1^{-}(\widetilde{K})$.

\begin{corollary}\colabel{36} Suppose that $q$ is not a root of unity in $k$. If $\ch k\neq 2$, then the set $$
\Bigl\{L_z(n,+), L_z(n,-)\Bigm|0\neq z\in k, n\in\mathbb{N}\cup\{0\}\Bigr\}.
$$
is a complete set of pairwise inequivalent finite-dimensional simple
weight $\wtu_q$-modules; while if $\ch k=2$, then the set $$
\Bigl\{L_z(n,+)=L_z(n,-)\Bigm|0\neq z\in k,
n\in\mathbb{N}\cup\{0\}\Bigr\}.
$$
is a complete set of pairwise inequivalent finite-dimensional simple
weight $\wtu_q$-modules.
\end{corollary}

\begin{proof} This follows from \leref{31} and \cite[Theorem 2.6]{Ja}.
\end{proof}

Radford in \cite{R3} constructed a class of simple Yetter--Drinfel'd
(shortly YD) modules for a graded Hopf algebra
$H=\oplus_{n=0}^{\infty}H_n$ with $H_0$ both commutative and
cocommutative. When $H$ is finitely graded over an algebraically
closed field and $H_0$ is the group algebra of a finite abelian
group, then all simple YD $H$-modules are in Radford's class of
simple YD modules. The Borel subalgebra ${\rm U}_q^{\leq 0}$,
denoted $H_\omega$ ($\omega=q^{-2}, a=K^{-1}$) in \cite{R3}, is a
simple pointed graded Hopf algebra, but not finitely graded. Thus we
don't know whether Radford's class of simple YD modules of ${\rm
U}_q^{\leq 0}$ \cite[Proposition 4,(b)]{R3}, parameterized by
$k^{\times}\times \mathbb{Z}$, contains all YD simple ${\rm
U}_q^{\leq 0}$-modules. However, we know that Radford's class forms
a proper subset of simple $D_q$-modules. Recall the Hopf algebra map
$\theta'$ we introduced in \leref{11}. Using $\theta'$ and noting
that the multiplication rule for our quantum double is compatible
with the multiplication rule given in \cite[Chapter IX, (4.3)]{K}
for the Drinfel'd quantum double of finite dimensional Hopf
algebras, one sees easily that every YD ${\rm U}_q^{\leq 0}$-module
naturally becomes a $D_q$-module. To keep in accordance with the
notations used in \cite[Proposition 4]{R3}, we set
$a=\widetilde{K}^{-1}, x=F$, $g=a^{l}, \omega=q^{-2}$, and let
$\beta: {\rm U}_q^{\leq 0}\rightarrow K$ be an algebra homomorphism.

\begin{proposition} With the notations as above and suppose that $q$ is not a root of unity in $k$, then

1) if $\beta(a)\neq\omega^{l+n}$ for any integer $n\geq 0$, then
$\lam^2\neq q^{2n}$ for any integer $n\geq 0$, where
$\lam:=\sqrt{\beta(a)}^{-1}q^{-l}, z:=\beta(a)q^{-2l}$, in this
case, the module $H_{\beta,kg}$ defined in \cite[Corollary 1]{R3} is
isomorphic (as $D_q$-module) to the infinite dimensional simple
$D_q$-module $M_{z}^{+}(\lam)$;

2) if  $\beta(a)=\omega^{l+n}$ for some integer $n\geq 0$, we set
$z=q^{-2(n+2l)}$, then the module $H_{\beta,kg}$ defined in
\cite[Corollary 1]{R3} is isomorphic (as $D_q$-module) to the
$(n+1)$-dimensional simple $D_q$-module $L_{z}(n,+)$ if
$q^{-n-2l}=\sqrt{z}$; or to the $(n+1)$-dimensional simple
$D_q$-module $L_{z}(n,-)$ if $q^{-n-2l}=-\sqrt{z}$.
\end{proposition}

\begin{proof} 1) Let $\lam:=\sqrt{\beta(a)}^{-1}q^{-l}, z:=\beta(a)q^{-2l}$. Then it is obvious that
$\beta(a)\neq\omega^{l+n}$ for any integer $n\geq 0$ if and only if
$\lam^2\neq q^{2n}$ for any integer $n\geq 0$. In this case, we know
that (by \coref{35}) $M_{z}^{+}(\lam)$ is a simple $D_q$-module. By
\cite[Proposition 4, (a)]{R3} and the formula given in the paragraph
below \cite[Proposition 4]{R3}, we have that $$
\left\{\begin{aligned}
\widetilde{K}^{-1}\bullet_{\beta}g&=a\bullet_{\beta}g=\beta(a)g,\,\,
K\bullet_{\beta}g=\varphi(K,K^{-l})g=q^{-2l}g,\\
E\bullet_{\beta}g&=\varphi(E,K^{-l})g=0.\end{aligned}\right.
$$
On the other hand, by the formula given in the paragraph above \coref{35}, we have that $$\left\{\begin{aligned}
\widetilde{K}^{-1}m_0&=z\lam^{-1}m_0=\beta(a)m_0,\,\,
K m_0=\lam m_0=q^{-2l}m_0,\\
Em_0&=0. \end{aligned}\right.
$$
By the universal property of the $D_q$-module $M_{z}^{+}(\lam)$ (see
\coref{35}), we deduce that the map which sends $m_0$ to $g$ can be
uniquely extended to a homomorphism $\eta$ from $M_{z}^{+}(\lam)$ to
$H_{\beta,kg}$. Comparing the action of $F$ on the basis
$\{x^i\bullet_{\beta}g\}_{i=0}^{\infty}$ given in the paragraph
below \cite[Proposition 4]{R3} and the action of $F$ on the basis
$\{m_i\}_{i=0}^{\infty}$ given in the paragraph above \coref{35}, we
know that $\eta(m_i)=x^i\bullet_{\beta}g$ for each $i\geq 0$, hence
$\eta$ is an isomorphism, as required.

2) We consider only the case where $q^{-n-2l}=\sqrt{z}$, the other
case is similar.

By the formula given in the paragraph below \cite[Proposition 4]{R3}, we have that $$
\left\{\begin{aligned}
\widetilde{K}^{-1}\bullet_{\beta}g&=a\bullet_{\beta}g=\omega^{l+n}g=q^{-2(l+n)}g,\\
K\bullet_{\beta}g&=\varphi(K,K^{-l})g=q^{-2l}g,\\
E\bullet_{\beta}g&=\varphi(E,K^{-l})g=0.\end{aligned}\right.
$$
On the other hand, by the formula given in the second paragraph below \coref{35}, we have that $$\begin{aligned}
\widetilde{K}^{-1}m_0&=q^{-2(l+n)},\,\,
K m_0=q^{-2l}m_0,\,\,
Em_0=0. \end{aligned}
$$
By the universal property of the $D_q$-modules $M_{z}^{+}(q^{-2l}),
L_{z}(n,+)$ (see \coref{35}), we deduce that the map which sends
$m_0$ to $g$ can be uniquely extended to a homomorphism $\eta'$ from
$M_{z}^{+}(q^{-2l})$ to $H_{\beta,kg}$ and hence gives rises to a
homomorphism $\eta'$ from $L_{z}(n,+)$ to $H_{\beta,kg}$. Comparing
the action of $F$ on the basis $\{x^i\bullet_{\beta}g\}_{i=0}^{n}$
given in the paragraph below \cite[Proposition 4]{R3} and the action
of $F$ on the basis $\{m_i\}_{i=0}^{n}$ given in the second
paragraph below \coref{35}, we know that
$\eta'(m_i)=x^i\bullet_{\beta}g$ for each $i\geq 0$, hence $\eta'$
is an isomorphism, as required. \end{proof}

It would be interesting to know if every simple $D_q$-module comes from a simple YD ${\rm U}_q^{\leq 0}$-module  when $q$
is not a root of unity. If this is the case, we would know all simple YD ${\rm U}_q^{\leq 0}$-modules. Taking the
advantage of the well established representation theory of the quantized enveloping algebra $U_q(\ksl_2)$, we can
easily obtain the decomposition of the tensor product of two finite dimensional $D_q$-modules while it might be
difficult for the YD module setting of Radford in \cite{R3}.

\begin{theorem} Suppose that $q$ is not a root of unity in $k$. Let $z,z'\in k^{\times}$.

Let $m,n\in\mathbb{N}\cup\{0\}$. Then there is a decomposition of
$\wtu_q$-modules: $$\begin{aligned} &L_z(m,\pm)\otimes
L_{z'}(n,\pm)\cong\begin{cases}
\oplus_{i=0}^{\min(m,n)}L_{zz'}(m+n-2i,+) &\text{if $\sqrt{z}\sqrt{z'}=\sqrt{zz'}$}\\
\oplus_{i=0}^{\min(m,n)}L_{zz'}(m+n-2i,-) &\text{if
$\sqrt{z}\sqrt{z'}=-\sqrt{zz'}$}
\end{cases}
\\
&L_z(m,\pm)\otimes L_{z'}(n,\mp)\cong\begin{cases}
\oplus_{i=0}^{\min(m,n)}L_{zz'}(m+n-2i,-)
&\text{if $\sqrt{z}\sqrt{z'}=\sqrt{zz'}$}\\
\oplus_{i=0}^{\min(m,n)}L_{zz'}(m+n-2i,+) &\text{if
$\sqrt{z}\sqrt{z'}=-\sqrt{zz'}$}
\end{cases}
\end{aligned}$$
\end{theorem}

\begin{proof} The theorem follows easily from \leref{tensor} and the
Clebsch--Gordan formula for $L(n,\pm)\otimes L(m,\pm)$.
\end{proof}

Many results for the quantized enveloping algebra $\ru_q$ have
their analogues for the algebra $\wtu_q$. For example, it is not
hard to show that the element
$C=FE+\DDF{Kq+\widetilde{K}^{-1}q^{-1}}{(q-q^{-1})^2}$ is equal to
$EF+\DDF{Kq^{-1}+\widetilde{K}^{-1}q}{(q-q^{-1})^2}$, and $C$ is
in the center of $\wtu_q$ (compare with \cite[(2.7)]{Ja}). In
\cite[\S2.13]{Ja}, the classification of the finite
dimensional simple $\ru_q$-modules are given when $q$ is a
primitive $l$-th root of unity with $l$ odd (the even case is also
similar). As a consequence, we have an analogous classification
result for the algebra $\wtu_q$. For example, when $q^2$ is a
primitive $d$-th root of unity in $k$ with $d>1$, we still have
the well-defined simple $\wtu_q$-modules $L_z(n,\pm)$  for $0\neq
z\in k$ and any integer $n$ with $0\leq n<d$.

\begin{lemma} \lelabel{39} Let $q^2$ be a primitive $d$-th root of unity in $k$ with $d>1$. If $M$ is a finite dimensional
simple weight $\wtu_q$-module such that both $E^{d}$ and $F^{d}$ act
as $0$ on $M$, then $M$ is isomorphic to one of the following
modules:
$$ Z_{0,z}^{\pm}(\lam),\,\,L_z(n,+),\,\,L_z(n,-),\,\,0\neq z\in
k,\,0\leq n< d,
$$
where $0\neq z\in k, 0\neq\lam\in k$ with $\lam^{2d}\neq 1$, and
$$Z_{0,z}^{\pm}(\lam):=M_{z}^{\pm}(\lam)/(\wtu_qm_d).
$$\end{lemma}

\bigskip

\section{Connections with the Drinfel'd double of the Taft algebra}

Throughout this section, we assume that $k$ is an algebraically closed field, $1<d\in\mathbb{N}$ and that $q^2\in k$
is a primitive $d$-th root of unity.

We consider the quantized enveloping algebra ${\ru}_q$. It is
well-known that the elements $E^d, F^d, K^d$ are central in the
algebra ${\ru}_q$. Let $\overline{\ru}_q^{\geq 0}$ (resp.
$\overline{\ru}_q^{\leq 0}$) be the quotient of the algebra
${\ru}_q^{\geq 0}$ (resp. ${\ru}_q^{\leq 0}$) modulo the ideal
generated by $E^d, K^d-1$ (resp. by $F^d, K^d-1$). It is well-known
that the ideal generated by $E^d, K^d-1$ (resp. by $F^d, K^d-1$) is
a Hopf ideal. Hence the algebra $\overline{\ru}_q^{\geq 0}$ (resp.
the algebra $\overline{\ru}_q^{\leq 0}$) is a quotient Hopf algebra
of ${\ru}_q^{\geq 0}$ (resp. of ${\ru}_q^{\leq 0}$). Recall the skew
Hopf pairing between ${\ru}_q^{\geq 0}$ and ${\ru}_q^{\leq 0}$
defined in Section 1.

\begin{lemma} \lelabel{41} The elements $E^d,K^d-1\in {\ru}_q^{\geq 0}, F^d,K^d-1\in {\ru}_q^{\leq 0}$ lie in the radical
of the skew Hopf pairing. Moreover, the induced skew Hopf pairing between $\overline{\ru}_q^{\geq 0}$ and
$\overline{\ru}_q^{\leq 0}$ is non-degenerate.
\end{lemma}

\begin{proof} For convenience, we still denote by $E^aK^b$ the canonical image of $E^aK^b\in{\ru}_q^{\geq 0}$ in
$\overline{\ru}_q^{\geq 0}$, and do the same for the elements $F^aK^b\in{\ru}_q^{\leq 0}$. Note that the monomials
$\bigl\{E^aK^b\bigr\}_{0\leq a,b<d}$ (resp. $\bigl\{F^aK^b\bigr\}_{0\leq a,b<d}$)
form a $k$-basis of $\overline{\ru}_q^{\geq 0}$ (resp.
$\overline{\ru}_q^{\leq 0}$). Recall that for the skew Hopf
pairing $\varphi$ between ${\ru}_q^{\geq 0}$ and ${\ru}_q^{\geq
0}$,
$$
\text{$\varphi(E^aK^b, F^{a'}K^{b'})=0$ unless $a=a'$ (cf. \cite[Proposition 1.2.3(d)]{L5})}
$$
With this in mind, the first statement of the lemma follows from a
direct verification. It remains to show that the induced skew Hopf
pairing is non-degenerate.

Suppose that $x:=\sum_{0\leq a,b<d}\lam_{a,b}E^aK^b\in\overline{\ru}_q^{\geq 0}$ (where $\lam_{a,b}\in k$ for each $a,b$)
lies in the radical of the induced skew Hopf pairing between $\overline{\ru}_q^{\geq 0}$ and $\overline{\ru}_q^{\leq 0}$.
We want to show that $\lam_{a,b}=0$ for all $a,b$.

Let $0\leq a<d$ be a fixed integer. By assumption, we have that $$
0=\varphi(x,F^{a}K^{b'})=\sum_{0\leq b<d}\lam_{a,b}\varphi(E^aK^b, F^aK^{b'}),\,\,\text{for $b'=0,1,2,\cdots,d-1$}.
$$
It is not hard to calculate that
$$
\varphi(E^aK^b, F^aK^{b'})=q^{-2bb'}[a]_{q}^{!}\Bigl(\frac{1}{1-q^2}\Bigr)^a.
$$
Since $q^2$ is a primitive $d$th root of unity, it follows that $[a]_{q}^{!}\neq 0$. Hence we get that $$
\sum_{0\leq b<d}\lam_{a,b}q^{-2bb'}=0,\,\,\,\,\text{for $b'=0,1,2,\cdots,d-1$}.
$$
Note that the coefficient matrix of the above system of linear equations is the Vandermonde matrix:
$$
\begin{pmatrix}
1& 1& 1& \cdots & 1\\
1& q^{-2}& q^{-4}& \cdots &q^{-2(d-1)}\\
1& (q^{-2})^2& (q^{-4})^{2}& \cdots &(q^{-2(d-1)})^{2}\\
\vdots & \vdots & \vdots & &\vdots \\
1& (q^{-2})^{d-1}& (q^{-4})^{d-1}&  \cdots &(q^{-2(d-1)})^{d-1}\\
\end{pmatrix},
$$
which has the non-zero determinant. It follows that $\lam_{a,b}=0$ for all
$0\leq a,b<d$ and hence $x=0$ as desired. In a similar way, one can prove that if
$y:=\sum_{0\leq a,b<d}\lam_{a,b}F^aK^b\in\overline{\ru}_q^{\leq 0}$ (where $\lam_{a,b}\in k$ for each $a,b$) lies in
the radical of the induced skew Hopf pairing between $\overline{\ru}_q^{\geq 0}$ and $\overline{\ru}_q^{\leq 0}$,
then $y=0$. This completes the proof of the lemma.
\end{proof}

Since $\overline{\ru}_q^{\leq 0}$ is of finite dimension, we have the following consequence of \leref{41}.

\begin{corollary} With the above induced skew Hopf pairing, the associated quantum double of $\overline{\ru}_q^{\geq 0}$
and $\overline{\ru}_q^{\leq 0}$ is isomorphic to the usual Drinfel'd double {\rm (cf. \cite{D2}, \cite{K})} of
$\overline{\ru}_q^{\leq 0}$ as a finite-dimensional $k$-Hopf algebra.
\end{corollary}

Denote by $\overline{\wtu}_q$ the quantum double of
$\overline{\ru}_q^{\geq 0}$ and $\overline{\ru}_q^{\leq 0}$ under
the above skew Hopf pairing. Note that the ideal generated by
$E^{d},\,\,F^d,\,\,K^d-1,\,\,\widetilde{K}^d-1$ is a Hopf ideal of
$\wtu_q$.

\begin{theorem} \thlabel{43} As a Hopf algebra, $\overline{\wtu}_q$ is isomorphic to the quotient of $\wtu_q$ modulo the
ideal generated by
$E^{d},\,\,F^d,\,\,K^d-1,\,\,\widetilde{K}^d-1$.
\end{theorem}

Note that $\overline{\ru}_q^{\geq 0}$ and $\overline{\ru}_q^{\leq 0}$ are isomorphic as $k$-Hopf algebras. Thus
$\overline{\ru}_q^{\leq 0}$ is a self-dual Hopf algebra. This Hopf algebra is usually called the Taft algebra,
denoted by $T_{d}(q^{-2})$, as it was constructed in \cite{T} as an interesting class of pointed Hopf algebras.

In \cite{C1}, Chen classified the irreducible representations of the
Drinfel'd double of $T_{d}(q^{-2})$ and studied their tensor
products. We remark that most of the results obtained in \cite{C1}
can be recovered easily from our \leref{tensor} and the discussion
below Lemma 3.3, and the corresponding known results for
${\ru_q}(\ksl_2)$. For example, our \leref{39} recovers the
classification of simple $D(T_{d}(q^{-2}))$-modules obtained in
\cite[Proposition 2.4, Theorem 2.5]{C1}. The decomposition formula
for the tensor product of two simple $D(T_{d}(q^{-2}))$-modules
obtained in \cite[Theorem 3.1]{C1} follows easily from our
\leref{tensor} and \cite[Theorem 4.5]{S} (see also \cite{P}).
Moreover, some results about finite dimensional indecomposable
representations of $D(T_d(q^{-2}))$ in \cite{C2} can be recovered
from \thref{43} and the results from \cite{X2}.

\section{Generalization to the case of arbitrary Cartan matrix}

Our main results in Section 3 allow a direct generalization to the
case of arbitrary Cartan matrix. To be precise, let
$A=(a_{i,j})_{1\leq i,j\leq n}$ be a $n\times n$ matrix with
entries in $\{-3,-2,-1,0,2\}$, $a_{i,i}=2$ and $a_{i,j}\leq 0$ for
$i\neq j$. Suppose $(d_1,\ldots,d_n)$ is a vector with entries
$d_i\in\{1,2,3\}$ such that the matrix $(d_ia_{i,j})$ is symmetric
and positive definite. Then $A$ is a Cartan matrix. Let
$\alpha_1,\ldots,\alpha_n$ be the set of simple roots in the
corresponding root system.

Let $k$ be a field. Let $q$ be an invertible element in $k$
satisfying $q^{2d_i}\neq 1$ for every $1\leq i\leq n$. The quantized enveloping algebra ${\ru}_q$ associated
to the Cartan matrix $A=(a_{i,j})_{1\leq i,j\leq n}$ (cf. \cite{D1}, \cite{J1} and \cite{J2})
is the associative $k$-algebra with generators $E_i, F_i, K_i, K_i^{-1} (1\leq i\leq n)$ and the
relations: $$
\begin{aligned}
& K_iK_j=K_jK_i,\quad K_iK_i^{-1}=1=K_i^{-1}K_i,\\
& K_iE_j=q^{d_ia_{i,j}}E_jK_i,\,\,K_iF_j=q^{-d_ia_{i,j}}F_jK_i,\\
& E_iF_j-F_jE_i=\delta_{i,j}\frac{K_i-K_i^{-1}}{q^{d_i}-q^{-d_i}},\\
& \sum_{r+s=1-a_{i,j}}(-1)^s\begin{bmatrix}1-a_{i,j}\\
s\end{bmatrix}_{q^{d_i}}E_i^rE_jE_i^s=0,\quad\text{if $i\neq j$},\\
& \sum_{r+s=1-a_{i,j}}(-1)^s\begin{bmatrix}1-a_{i,j}\\
s\end{bmatrix}_{q^{d_i}}F_i^rF_jF_i^s=0,\quad\text{if $i\neq j$}.
\end{aligned}
$$

${\ru}_q$ is a Hopf algebra with comultiplication, counit and antipode given by:
$$\begin{aligned} &\Delta(E_i)=E_i\otimes 1+K_i\otimes
E_i,\,\Delta(F_i)=F_i\otimes K_i^{-1}+1\otimes
F_i,\,\Delta(K_i)=K_i\otimes K_i, \\
& \varepsilon(E_i)=0=\varepsilon(F_i),\,\,\varepsilon(K_i)=1=\varepsilon(K_i^{-1}),\\
& S(E_i)=-K_i^{-1}E_i,\,\,S(F_i)=-F_iK_i,\,\,S(K_i)=K_i^{-1}.
\end{aligned}
$$
Let ${\ru}_q^{+}$ (resp. ${\ru}_q^{-}$) be the $k$-subalgebra of
${\ru}_q$ generated by $E_i, 1\leq i\leq n$ (resp. by $F_i, 1\leq
i\leq n$). Let ${\ru}_q^{0}$ be the $k$-subalgebra of ${\ru}_q$
generated by $K_i, K_i^{-1}, 1\leq i\leq n$. Let ${\ru}_q^{\geq
0}:={\ru}_q^{+}{\ru}_q^{0}$, ${\ru}_q^{\leq
0}:={\ru}_q^{-}{\ru}_q^{0}$. Then both ${\ru}_q^{\geq 0}$ and
${\ru}_q^{\leq 0}$ are Hopf $k$-subalgebras of ${\ru}_q$. For any
monomials $$(\prod_{1\leq i\leq n}E_i^{a_i})(\prod_{1\leq i\leq
n}K_i^{b_i})\in {\ru}_q^{+},\quad (\prod_{1\leq i\leq
n}F_i^{a_i})(\prod_{1\leq i\leq n}K_i^{b_i})\in{\ru}_q^{-}, $$ we
endow them the weights $\sum_{i=1}^{n}a_i\alpha_i,
-\sum_{i=1}^{n}a_i\alpha_i$ respectively. Like the Hopf pair
$(U_q(\ksl_2)^{\geq 0}),U_q(\ksl_2)^{\leq 0})$,  there exists a
unique pairing $\varphi: {\ru}_q^{\geq 0}\times {\ru}_q^{\leq
0}\rightarrow k$ (see \cite{L5}, \cite{Jo} and \cite{X}) such that
$$\begin{aligned} (1)
&\,\,\,\varphi(1,1)=1,\,\,\varphi(1,K_i)=1=\varphi(K_i,1),\,\text{for $1\leq i\leq n$,}\\
(2) &\,\,\,\varphi(x,y)=0,\,\,\text{if $x,y$ are homogeneous with different weights,}\\
(3) &\,\,\,\varphi(E_i,F_j)=\delta_{i,j}\frac{1}{q^{2d_i}-1},\,\text{for $1\leq i,j\leq n$,}\\
(4) &\,\,\,\varphi(K_i,K_j)=q^{d_ia_{i,j}},\,\,
\varphi(K_i,K_j^{-1})=q^{-d_ia_{i,j}},\,\text{for $1\leq i,j\leq n$,}\\
(5) &\,\,\,\varphi(x,y'y'')=\varphi(\Delta^{\op}(x),y'\otimes
y''),\,\,\text{for all $x\in{\ru}_q^{\geq 0},\,
y',y''\in{\ru}_q^{\leq 0}$},\\
(6) &\,\,\,\varphi(xx',y'')=\varphi(x\otimes
x',\Delta(y'')),\,\,\text{for all $x,x'\in{\ru}_q^{\geq 0},\,
y''\in{\ru}_q^{\leq 0}$},\\
(7) &\,\,\,\varphi(S(x),y)=\varphi(x,S^{-1}(y)),\,\,\text{for all
$x\in{\ru}_q^{\geq 0},\, y\in{\ru}_q^{\leq 0}$}.
\end{aligned}
$$
In other words, $({\ru}_q^{\geq 0}, {\ru}_q^{\leq 0}, \varphi)$
forms a skew Hopf pairing. Thus (as in Section 1) we can make
$D({\ru}_q^{\geq 0}, {\ru}_q^{\leq 0}):={\ru}_q^{\geq 0}\otimes
{\ru}_q^{\leq 0}$ into a Hopf $k$-algebra,  called the quantum
double of $({\ru}_q^{\geq 0}, {\ru}_q^{\leq 0},\varphi)$. For
simplicity, we write $\wtu_q$ instead of $D({\ru}_q^{\geq 0},
{\ru}_q^{\leq 0})$.

\begin{theorem} As a $k$-algebra, $\wtu_q$ can be presented by the generators $$E_i, F_i, K_i, K_i^{-1},
\widetilde{K}_i, \widetilde{K}_i^{-1},\,\, (1\leq i\leq n), $$ and the following relations:
$$\begin{aligned}
& K_iK_j=K_jK_i,\quad K_iK_i^{-1}=1=K_i^{-1}K_i,\\
& \widetilde{K}_i\widetilde{K}_j=\widetilde{K}_j\widetilde{K}_i,\quad \widetilde{K}_i\widetilde{K}_i^{-1}=1=
\widetilde{K}_i^{-1}\widetilde{K}_i,\quad K_i\widetilde{K}_j=\widetilde{K}_jK_i,\\
& K_iE_j=q^{d_ia_{i,j}}E_jK_i,\,\,K_iF_j=q^{-d_ia_{i,j}}F_jK_i,\\
& \widetilde{K}_iE_j=q^{d_ia_{i,j}}E_j\widetilde{K}_i,\,\,\widetilde{K}_iF_j=q^{-d_ia_{i,j}}F_j\widetilde{K}_i,\\
& E_iF_j-F_jE_i=\delta_{i,j}\frac{K_i-\widetilde{K}_i^{-1}}{q^{d_i}-q^{-d_i}},\end{aligned} $$

$$\begin{aligned}
& \sum_{r+s=1-a_{i,j}}(-1)^s\begin{bmatrix}1-a_{i,j}\\
s\end{bmatrix}_{q^{d_i}}E_i^rE_jE_i^s=0,\quad\text{if $i\neq j$},\\
& \sum_{r+s=1-a_{i,j}}(-1)^s\begin{bmatrix}1-a_{i,j}\\
s\end{bmatrix}_{q^{d_i}}F_i^rF_jF_i^s=0,\quad\text{if $i\neq j$}.
\end{aligned}
$$
\end{theorem}

Let $M$ be a $\wtu_q$-module such that $\End_{\wtu_q}(M)=k$. Note
that the elements $K_i\widetilde{K}_i^{-1}, i=1,2,\cdots,n$, are
invertible central elements in $\wtu_q$. Therefore, there is a
vector $\vec{z}=(z_1,\cdots,z_n)\in (k^{\times})^n$, such that for
every $1\leq i\leq n$, $K_i\widetilde{K}_i^{-1}$ acts as the scalar
$z_i$ on $M$. For each $0\neq z\in k$, we fix a square root
${z}^{1/2}$ of $z$. Let $\pi_{\vec z}^{+}$ be the $k$-algebra
homomorphism $\wtu_q\rightarrow{\ru}_q$ which is defined on
generators by $$ \pi_{\vec z}^{+}(E_i)=z_i^{1/2}E_i,\,\,\pi_{\vec
z}^{+}(F_i)=F_i,\,\,\pi_{\vec z}^{+}(K_i)=z_i^{1/2}K_i,\,\,
\pi_{\vec z}^{+}(\widetilde{K}_i)=z_i^{-1/2}K_i,
$$
for every $1\leq i\leq n$. It is easy to check that $\pi_{\vec
z}^{+}$ is well-defined. Moreover, the kernel of $\pi_{\vec z}^{+}$,
which is the ideal generated by $K_i\widetilde{K}_i^{-1}-z_i,
i=1,2,\cdots,n$, annihilates the module $M$. It follows that $M$
naturally becomes a module over the algebra ${\ru}_q$. Note that
$\pi_{\vec z}^{+}$ is in general not a Hopf algebra map unless
${\vec z}=(1,1,\cdots,1)$.

We call a $\wtu_q$-module $M$ a weight $\wtu_q$-module if
$K_1,\cdots,K_n,\widetilde{K}_1,\cdots,\widetilde{K}_n$ all act
semisimply on $M$. In that case, each $K_{i}\widetilde{K}_{i}^{-1}$
acts semisimply on $M$ as well. Similarly, we call a
${\ru}_q$-module $N$ a weight ${\ru}_q$-module if $K_1, K_2, \cdots,
K_n$ all act semisimply on $N$.

\begin{lemma} Every finite dimensional simple
(resp. indecomposable weight) $\wtu_q$-module is the pull-back of a
finite dimensional simple (resp. indecomposable weight)
${\ru}_q$-module through the algebra homomorphism $\pi_{\vec z}^{+}$
for some $\vec{z}=(z_1,\cdots,z_n)\in (k^{\times})^n$.
\end{lemma}

Let $M$ be a ${\ru}_q$-module. Let $\vec{z}=(z_1,\cdots,z_n)\in
(k^{\times})^n$. We use $M_{\vec z}^{+}$ to denote the pull-back of
$M$ through the algebra homomorphism $\pi_{\vec z}^{+}$. Let
$\varepsilon_{\vec z}^{+}$ be the one-dimensional representation of
$\wtu_q$ which is defined on generators by $$ \varepsilon_{\vec
z}^{+}(E_i)=0=\varepsilon_{\vec z}^{+}(F_i),\,\,\varepsilon_{\vec
z}^{+}(K_i)=z_i^{1/2},\,\, \varepsilon_{\vec
z}^{+}(\widetilde{K}_i)=z_i^{-1/2},\,i=1,2,\cdots,n.
$$
It is easy to check that $\varepsilon_{\vec z}^{+}$ is well-defined.

\begin{theorem} \thlabel{62} The category $\widetilde{\mathcal{C}}$ of finite dimensional weight $\wtu_q$-modules is
equivalent to a direct sum of $|(k^{\times})^n|$ copies of the
category $\mathcal{C}$ of finite dimensional weight
${\ru}_q$-modules.
\end{theorem}

With those one dimensional representations $\varepsilon_{\vec
z}^{+}$ in mind, one can also formulate a version of \leref{tensor}
in the context of arbitrary Cartan matrix. As before, this provides
an easy solution to the problem of decomposing the tensor product of
certain $\wtu_q$-modules, i.e., reducing them to the corresponding
problem for ${\ru}_q$-modules.

\bigskip\bigskip
\centerline{ACKNOWLEDGMENTS} \hskip\parindent The first author would
like to thank the School of Mathematics, Statistics and Computer
Science, Victoria University of Wellington for their hospitality
during his visit in 2005. He is grateful to the URF of VUW  and the
Program NCET as well as NSFC (Project 10401005) for the financial
support. The second author is supported by the Marsden Fund.

\end{document}